\documentclass[12pt]{amsart}
\usepackage{amsmath}
\usepackage{amsfonts}
\usepackage{amsthm}
\usepackage{amssymb}
\usepackage{tikz-cd}
\usepackage{graphics}
\usepackage{array}
\usepackage{color}
\usepackage{wasysym}
\usepackage{hyperref}
\usepackage{pdfsync}
\usepackage{mathrsfs}
\usepackage{mathtools}
\usepackage{needspace}
\usepackage[alphabetic]{amsrefs}
\usepackage{geometry}
\usepackage{bbm}
\usepackage[all,cmtip]{xy}
\usepackage{tikz}\usetikzlibrary{intersections}
\usepackage{extpfeil}

\newcommand{\R}{{\mathbb R}}
\newcommand{\N}{{\mathbb N}}

\newcommand{\C}{{\mathbb C}}
\newcommand{\Q}{{\mathbb Q}}

\DeclareMathOperator{\Image}{Im}
\DeclareMathOperator{\Pic}{Pic}

\DeclareMathOperator{\Coker}{Coker}
\DeclareMathOperator{\codim}{codim}
\DeclareMathOperator{\Supp}{Supp}

\DeclareMathOperator{\GV}{GV}
\DeclareMathOperator{\IT}{IT_0}

\DeclareMathOperator{\id}{id}

\DeclareMathOperator{\tors}{tors}
\DeclareMathOperator{\Sym}{Sym}

\theoremstyle{plain}

\newtheorem{thm}{Theorem}[section]
\newtheorem{conj}[thm]{Conjecture}
\newtheorem{prop}[thm]{Proposition}
\newtheorem{coro}[thm]{Corollary}
\newtheorem{lemma}[thm]{Lemma}

\theoremstyle{definition}

\newtheorem{defi}[thm]{Definition}

\newtheorem{rem}[thm]{Remark}
\newtheorem*{ac}{Acknowledgements}

\newcommand{\cO}{\mathcal{O}}

\newcommand{\cI}{\mathcal{I}}
\newcommand{\cM}{\mathcal{M}}
\newcommand{\cF}{\mathcal{F}}

\newcommand{\cE}{\mathcal{E}}

\newcommand{\cR}{\mathcal{R}}

\newcommand{\cK}{\mathcal{K}}
\newcommand{\cN}{\mathcal{N}}
\newcommand{\cQ}{\mathcal{Q}}
\newcommand{\cU}{\mathcal{U}}

\DeclareMathOperator{\Spec}{Spec}

\DeclareMathOperator{\Hom}{Hom}

\DeclareMathOperator{\Vol}{Vol}
\DeclareMathOperator{\ord}{ord}

\DeclareMathOperator{\coef}{coef}

\DeclareMathOperator{\Bl}{Bl}

\numberwithin{equation}{section}

\begin{document}

\title[Pushforwards of lc pairs]{Pushforwards of lc pairs under morphisms to abelian varieties}

\author{Fanjun Meng}
\address{Department of Mathematics, University of California San Diego, 9500 Gilman Drive \# 0112, La Jolla, CA 92093-0112, USA}
\email{f2meng@ucsd.edu}
    

\begin{abstract}
Let $f$ be a morphism from a log canonical pair $(X,\Delta)$ to an abelian variety $A$ and $D\sim_{\Q}m(K_X+\Delta)$ a Cartier divisor on $X$ with $m>1$ a rational number. We prove that $f_*\cO_X(D)$ is globally generated after pullback by an isogeny. Our methods involve birational geometry, mixed Hodge modules and singular Hermitian metrics with the minimal extension property. Moreover, we provide simple counterexamples where the above result fails for plt pairs when $m=1$ and where the existence of Chen--Jiang decompositions fails for plt pairs for any $m\ge1$.
\end{abstract}

\maketitle
\setcounter{tocdepth}{1}
\tableofcontents

\section{Introduction}\label{sec:introduction}

In this paper, we study pushforwards of log canonical pairs (i.e. lc pairs) under morphisms to abelian varieties over $\C$. When $X$ is a smooth projective variety, positivity properties of pushforwards of pluricanonical bundles under morphisms to abelian varieties have been studied extensively. For example, they are GV-sheaves, their cohomological support loci are finite unions of torsion subvarieties, they are globally generated after pullback by an isogeny, and they have Chen--Jiang decompositions, see \cites{GL87,GL91,Sim93,Hac04,Lai11,PP11a,PS14b,PPS17,CJ18,LPS20}. When $(X,\Delta)$ is klt, the corresponding results regarding global generation after pullback by an isogeny and the existence of Chen--Jiang decompositions are independently proved by \cite{Men21} and \cite{Jia22}. Both authors in their papers ask whether similar statements hold for lc pairs. In this paper, we give a complete answer to this question. Interestingly, the global generation statement holds when $m>1$ but does not hold when $m=1$. The existence of Chen--Jiang decompositions does not hold for any $m\ge1$.

For the clarity of our discussion, we record the following conjectures here. They are the main topics studied in this paper. It is natural to consider them, especially after the work \cite{LPS20} for smooth projective varieties. We claim no originality or priority of proposing these conjectures.

\begin{conj}\label{conj:global generation}
Let $f$ be a morphism from an lc pair $(X, \Delta)$ to an abelian variety $A$, $m\geq1$ a rational number and $D$ a Cartier divisor on $X$ such that $D\sim_{\Q}m(K_X+\Delta)$. Then there exists an isogeny $\varphi\colon A'\to A$ such that $\varphi^*f_*\mathcal{O}_X(D)$ is globally generated.	
\begin{center}
\begin{tikzcd}
			X' \arrow[r, "\varphi'"] \arrow[d, "f'"] & X \arrow[d, "f"] \\
			 A' \arrow[r, "\varphi" ] & A 
\end{tikzcd}
\end{center}
\end{conj}

An isogeny is a finite surjective homomorphism of abelian varieties. Thus $A'$ is an abelian variety. The morphism $f$ does not have to be surjective.

\begin{conj}\label{conj:CJ decomposition}
Let $f$ be a morphism from an lc pair $(X, \Delta)$ to an abelian variety $A$, $m\geq1$ a rational number and $D$ a Cartier divisor on $X$ such that $D\sim_{\Q}m(K_X+\Delta)$. Then $f_*\mathcal{O}_X(D)$ admits a finite direct sum decomposition
\[
f_*\mathcal{O}_X(D)\cong \bigoplus_{i\in I}(\alpha_i\otimes p_i^*\mathcal{F}_i),
\]
where each $A_i$ is an abelian variety, each $p_i\colon A\to A_i$ is a fibration, each $\mathcal{F}_i$ is a nonzero M-regular coherent sheaf on $A_i$, and each $\alpha_i\in\Pic^0(A)$ is a torsion line bundle.
\end{conj}

The above decomposition is called the \emph{Chen--Jiang decomposition} for $f_*\mathcal{O}_X(D)$, see Definition \ref{defi:CJ decomposition}. M-regularity proposed in \cite{PP03} is an analogue of ampleness for coherent sheaves on abelian varieties and stronger than ampleness in general. However, these two notions are equivalent for line bundles on abelian varieties. Since pullbacks of ample line bundles are semiample, one can think of Conjecture \ref{conj:CJ decomposition} as saying $f_*\mathcal{O}_X(D)$ is semiample in some sense.

\begin{conj}\label{conj:canonical direct summand}
Let $f$ be a morphism from an lc pair $(X, \Delta)$ to an abelian variety $A$, $m\geq1$ a rational number and $D$ a Cartier divisor on $X$ such that $D\sim_{\Q}m(K_X+\Delta)$. Then there exists a generically finite surjective morphism $h\colon Z\to X$ from a smooth projective variety $Z$ such that $f_*\mathcal{O}_X(D)$ is a direct summand of $(f\circ h)_*\mathcal{O}_Z(K_Z)$.
\end{conj}

It follows from the proof of \cite{LPS20}*{Theorem 5.1} that Conjecture \ref{conj:CJ decomposition} implies Conjecture \ref{conj:global generation}. Conjecture \ref{conj:canonical direct summand} implies Conjecture \ref{conj:CJ decomposition} by \cite{PPS17}*{Theorem A} and \cite{LPS20}*{Proposition 3.6}. In summary, we have
\[
\mathrm{Conjecture\,\ref{conj:canonical direct summand}}\Rightarrow\mathrm{Conjecture\,\ref{conj:CJ decomposition}}\Rightarrow\mathrm{Conjecture\,\ref{conj:global generation}}.
\]
For klt pairs, these three conjectures have been proved to be equivalent before they are proved to be true by \cite{Men21}*{Proposition 3.2}. The equivalence is a useful tool for \cite{Men21} to establish these conjectures for klt pairs, which would help us if we have an lc version of it. The proof of this equivalence depends on several techniques from birational geometry and the fact that these conjectures hold for klt pairs when $m=1$. However, for lc pairs, the situation is not clear even when $m=1$. In fact, all these three conjectures fail for $m=1$ and Conjectures \ref{conj:CJ decomposition} and \ref{conj:canonical direct summand} even fail for any $m>1$, which indicates that Conjecture \ref{conj:global generation} is the suitable one to study in general. Moreover, they fail in a very nice setting where the pairs considered are log smooth plt surface pairs by Theorems \ref{thm:counterexample m=1} and \ref{thm:counterexample m>=1}. Thus these three conjectures are not equivalent in general. Plt pairs are more general than klt pairs. Plt pairs are lc.

In this paper, we show that Conjecture \ref{conj:global generation} is true for any $m>1$ and is false for $m=1$, while Conjectures \ref{conj:CJ decomposition} and \ref{conj:canonical direct summand} are false for any $m\ge1$. During our studies of these three conjectures, it becomes evident to us that Conjecture \ref{conj:global generation} should be the first one to study and it is more promising to establish positive results regarding it. Before we come to the counterexamples, let us start with the more significant positive result for lc pairs.

\begin{thm}\label{thm:gg for lc}
Conjecture \ref{conj:global generation} holds for any rational number $m>1$.
\end{thm}

Theorem \ref{thm:gg for lc} implies that $f_*\mathcal{O}_X(D)$ has strong positivity, which has many applications in algebraic geometry. For example, it directly implies the subadditivity of Kodaira dimension of fibrations over varieties of maximal Albanese dimension for lc pairs. For klt pairs, the aforementioned subadditivity is proved by \cite{CP17}*{Theorem 1.1} (see also \cite{HPS18}). It also directly implies that $f_*\mathcal{O}_X(D)$ admits a singular Hermitian metric with semipositive curvature and the minimal extension property when $f$ is surjective and $m>1$ for lc pairs, see Corollary \ref{coro:existence of singular Hermitian metrics with semipositive curvature and the minimal extension property}. When $m=1$, $f_*\mathcal{O}_X(D)$ may not admit any such metric by Theorem \ref{thm:counterexample m=1}. For smooth varieties, the existence of such a metric on $f_*\mathcal{O}_X(D)$ is proved by \cite{HPS18}*{Theorem 27.1} for any $m\ge1$. For klt pairs, it is proved by \cite{Men21}*{Lemma 2.12 and Theorem 1.1} for any $m\ge1$.

As a direct corollary of Theorem \ref{thm:gg for lc}, Conjecture \ref{conj:global generation} also holds when $m>1$ is a real number, $(X, \Delta)$ is an lc pair with real coefficients, and $D$ is a Cartier divisor on $X$ such that $D\sim_{\R}m(K_X+\Delta)$, see Corollary \ref{coro:real coefficient gg for lc}.

Conjecture \ref{conj:global generation} for lc pairs is more difficult than it for smooth varieties or klt pairs due to the lack of suitable techniques and strategies. Many of the previous methods break in our lc setting. Let us think heuristically about the difficulties in proving Theorem \ref{thm:gg for lc}:

\begin{enumerate}
\item Contrary to some potential intuitions, this global generation problem becomes more difficult when going from $m=1$ to $m>1$ even in the smooth case. General expectations from birational geometry tell us that $mK_X$ has stronger and stronger positivity when we increase $m$, and thus positivity problems regarding $mK_X$ should become easier when $m$ is large. Similar expectations do not apply easily to the pushforwards of these objects. Naively, there are no simple connections between $f_*\cO_X(mK_X)$ and various kinds of powers of $f_*\cO_X(K_X)$. Practically, all the positivity results regarding the sheaves of pushforwards considered here are always proved first for the canonical case $m=1$ and then for the pluricanonical case $m>1$. Mathematically, Hodge theoretic methods are powerful enough to prove these results for the canonical case since pushforwards of canonical bundles appear naturally in Hodge theory, while pushforwards of pluricanonical bundles do not. This is one of the reasons why various different techniques outside Hodge theory are needed to prove the results for the pluricanonical case.

\item For the positivity results considered here, the pluricanonical case is always proved by reducing it to the canonical case via various techniques. For generic vanishing, the reduction is done in \cite{PS14b} by proving a Koll\'ar-type vanishing theorem in the pluricanonical case and using techniques from birational geometry. For the description of the cohomological support loci, the reduction is done in \cites{Shi16,HPS18} by techniques from birational geometry. For the global generation property and the existence of Chen--Jiang decompositions, the reduction is done in \cites{LPS20,Men21} by using singular Hermitian metrics with semipositive curvature and the minimal extension property, techniques from birational geometry and induction on the dimension of the base abelian variety. However, for our global generation problem in the lc setting, the reduction fails since the base case of the reduction when $m=1$ already fails by Theorem \ref{thm:counterexample m=1}. Moreover, the established results for $m=1$ are needed in several places during the reduction from $m>1$ to $m=1$. Thus the previous strategies to prove the case when $m>1$ by reducing it to the case when $m=1$ do not work in our scenario.
\end{enumerate}

\subsection*{Sketch of the proof of Theorem \ref{thm:gg for lc}} By the above discussion, we can see that new strategies and techniques are needed in our work. Let us discuss very briefly some rough ideas of our strategies and techniques. First, without loss of generality, we can assume that $(X, \Delta)$ is log smooth lc. We prove Theorem \ref{thm:gg for lc} by induction first on $\dim A$ and then on the number of lc centers of $(X, \Delta)$. If the number of lc centers of $(X, \Delta)$ is $0$, then $(X, \Delta)$ is klt. We conclude by \cite{Men21}*{Theorem 1.1}. Next, we assume that the number of lc centers of $(X, \Delta)$ is $\ge1$. Choose a minimal lc center $T$ of $(X, \Delta)$, and consider the pushforward of the natural restriction sequence 
\[
 0\longrightarrow f_*(\cI_T\otimes\cO_X(D))\longrightarrow f_*\cO_X(D)
 \longrightarrow(f|_T)_*\mathcal O_T(D|_T).
\]
Working with minimal lc centers gives us a lot of advantages in our proof. Since $T$ is minimal, $D|_T\sim_{\Q}m(K_T+\Delta_T)$ where $(T, \Delta_T)$ is a klt pair by subadjunction. Thus by \cite{Men21}*{Theorem 1.1}, $(f|_T)_*\mathcal O_T(D|_T)$ is globally generated after pullback by an isogeny. However, it is not very helpful. First, the sequence above is not exact. Most importantly, the global sections of $(f|_T)_*\mathcal O_T(D|_T)$ may not lift to global sections of $f_*\cO_X(D)$, which is a notorious obstacle in these kinds of problems on positivity. One of the key points we notice is that we should consider
\[
\mathcal Q:=\Image\bigl(f_*\mathcal O_X(D)\longrightarrow(f|_T)_*\mathcal O_T(D|_T)\bigr),
\]
which seems natural, but is proved to be a correct way to deal with the global generation problem in the lc setting. It turns out to be quite difficult to make this observation work.

Now we have the following exact sequence
\[
0\longrightarrow f_*(\cI_T\otimes\cO_X(D))\longrightarrow f_*\cO_X(D)
 \longrightarrow\cQ\longrightarrow 0.
\]
By the induction hypothesis regarding the number of lc centers, it is not hard to see that $f_*(\cI_T\otimes\cO_X(D))$ is globally generated after pullback by an isogeny. By Proposition \ref{prop:lifting from the lc locus} whose proof uses extension from the non-lc locus, we can lift the global sections of $\cQ$ to the global sections of $f_*\cO_X(D)$ when $m>1$. Thus by considering $\cQ$, we get both the exactness of the relevant sequence and the lifting of the relevant global sections. Moreover, we prove that $\cQ$ has the similar nice positivity properties as $(f|_T)_*\mathcal O_T(D|_T)$ does. If $m=1$, $\cQ$ is a canonical direct summand (see Definition \ref{defi:canonical direct summand}) and is globally generated after pullback by an isogeny by Proposition \ref{prop:gg for lc restricted image m=1}. Three important ingredients in the proof of this proposition are Viehweg's cyclic covering trick, functorial resolutions of singularities and the theory of mixed Hodge modules. If $m>1$, $\mathsf T(\cQ)$, a canonically defined subsheaf of $\cQ$ (see Definition \ref{defi:continuous evaluation morphism}), is a canonical direct summand and is globally generated after pullback by an isogeny by Proposition \ref{prop:gg for lc restricted image m>1}. Moreover, for any $m\ge1$, $\cQ$ admits a singular Hermitian metric with semipositive curvature and the minimal extension property and every nonzero morphism $\cQ\to\cO_A$ splits by Theorem \ref{thm:image splitting}. With all these results in hand, we can conclude the proof of Theorem \ref{thm:gg for lc} with additional work regarding positivity properties of coherent sheaves on abelian varieties.

We already know that Conjecture \ref{conj:global generation} fails when $m=1$ by Theorem \ref{thm:counterexample m=1}. Let us take a closer look at why our proof of Theorem \ref{thm:gg for lc} fails when $m=1$. Actually, all the parts of our proof work perfectly when $m=1$ except Proposition \ref{prop:lifting from the lc locus} regarding the extension of global sections. Its proof essentially uses $m>1$. This proposition holds when $m>1$ but fails easily when $m=1$. Here is a simple counterexample. Let $A$ be an elliptic curve with origin $o$, $X=A$, $f=\id_A$, $D=\Delta=o$ and $T=\{o\}$. Then $(X,\Delta)$ is log smooth plt, $T$ is its unique
minimal lc center, and $D\sim K_X+\Delta$. The pushforward of the natural restriction sequence is
\[
0\longrightarrow\cO_A\longrightarrow\cO_A(o)\longrightarrow\cO_A(o)|_T\longrightarrow0.
\]
In our notation, $f_*\cO_X(D)=\cO_A(o)$, $\cQ=\cO_A(o)|_T$ and $\mathsf T(\cO_A)=\cO_A$. The map
\[
 H^0(A,f_*\cO_X(D))\longrightarrow H^0(A,\cQ)
\]
is not surjective.

\bigskip

The following theorem regarding the existence of singular Hermitian metrics on $\cQ$ with semipositive curvature and the minimal extension property is one of the central ingredients in the proof of Theorem \ref{thm:gg for lc}.

\begin{thm}[=Theorem \ref{thm:image splitting}] \label{main thm:image splitting}
Let $(X,S+B)$ be a log smooth lc pair with $S=\lfloor S+B\rfloor$, $f:X\to A$ a morphism to an abelian variety $A$, $m\ge1$ a rational number and $D$ a Cartier divisor on $X$ such that $D\sim_{\Q}m(K_X+S+B)$. Let $T$ be a reduced union of minimal lc centers of the same dimension such that each of them dominates $A$. Then the sheaf
\[
\mathcal Q:=\Image\bigl(f_*\mathcal O_X(D)\longrightarrow(f|_T)_*\mathcal O_T(D|_T)\bigr)
\]
admits a singular Hermitian metric with semipositive curvature and the minimal extension property.
\end{thm}

The theory of singular Hermitian metrics with semipositive curvature and the minimal extension property developed in \cite{HPS18} and based on the Ohsawa-Takegoshi extension theorem with sharp estimates plays an important role here, see also \cites{BP08,CP17,PT18}. In the klt setting, the situation is nice and we know that $f_*\mathcal O_X(D)$ admits a singular Hermitian metric with semipositive curvature and the minimal extension property by \cite{Men21}*{Lemma 2.12 and Theorem 1.1} for any $m\ge1$ when $f$ is surjective. In the lc setting, we do not know the existence of such a metric on $f_*\mathcal O_X(D)$ a priori. In fact, when $m>1$, this existence is a direct corollary of Theorem \ref{thm:gg for lc} \emph{only after} we have proved Theorem \ref{thm:gg for lc}. When $m=1$, $f_*\mathcal O_X(D)$ may not admit any such metric by Theorem \ref{thm:counterexample m=1}. Thus instead, we prove that its restricted image $\cQ$ admits such a metric when each irreducible component of $T$ dominates $A$ for any $m\ge1$. On a complex manifold, consider a torsion-free quotient sheaf of a torsion-free coherent sheaf admitting a singular Hermitian metric with semipositive curvature and the minimal extension property. It is not hard to see by definition that the induced quotient metric on the quotient sheaf satisfies the same properties. In our scenario, $\cQ$ is a torsion-free quotient sheaf of $f_*\mathcal O_X(D)$. If we knew that $f_*\mathcal O_X(D)$ admits such a metric, Theorem \ref{main thm:image splitting} would follow. However, we do not know that a priori. Instead, we establish Theorem \ref{main thm:image splitting} without knowing the existence of such a metric on $f_*\mathcal O_X(D)$. Note that when $m=1$, an interesting thing is that $\cQ$ admits a singular Hermitian metric with semipositive curvature and the minimal extension property while $f_*\mathcal O_X(D)$ may not admit any such metric by Theorem \ref{thm:counterexample m=1}.

We emphasize that by the above discussion, knowing only Theorem \ref{main thm:image splitting} is not enough for us to deduce Theorem \ref{thm:gg for lc} since we need the extension of global sections and much more ingredients. In the end, when $m=1$, Theorem \ref{main thm:image splitting} holds but Conjecture \ref{conj:global generation} fails by Theorem \ref{thm:counterexample m=1}.

\begin{rem}
Theorem \ref{main thm:image splitting} and many other results in this paper hold in the more general setting when the base variety is an arbitrary smooth (not necessarily quasi-projective) variety and not necessarily an abelian variety by the same methods as in this paper. We can prove that they still hold for the restricted image $\cQ$ when the base variety is an arbitrary smooth variety $Y$ and $D\sim_{\Q}m(K_{X/Y}+S+B)$. We will discuss them in our future work.
\end{rem}

Next, we give the promised two simple counterexamples for these conjectures.

\begin{thm}\label{thm:counterexample m=1}
Conjecture \ref{conj:global generation} is false for $m=1$. The counterexample is a log smooth plt surface pair where $K_X+\Delta$ is nef, big, and semiample, and $f_*\mathcal{O}_X(D)$ does not admit any singular Hermitian metric with semipositive curvature and the minimal extension property.
\end{thm}

\begin{thm}\label{thm:counterexample m>=1}
Conjectures \ref{conj:CJ decomposition} and \ref{conj:canonical direct summand} are false for any rational number $m\ge1$. The counterexample is a log smooth plt surface pair where $K_X+\Delta$ is nef, big, and semiample, Conjecture \ref{conj:global generation} holds for $\pi_*\cO_X(D)$, and $\pi_*\cO_X(D)$ admits a singular Hermitian metric with semipositive curvature and the minimal extension property.
\end{thm}

By \cite{PPS17}*{Theorem A} and \cite{LPS20}*{Proposition 3.6}, it is not hard to see from Theorem \ref{thm:counterexample m>=1} that a weaker version of Conjecture \ref{conj:canonical direct summand} is still false for any $m\ge1$ even when we do not put any restriction on the morphism $h\colon Z\to X$ from a smooth projective variety $Z$ (i.e. we do not ask $h$ to be generically finite or surjective in Conjecture \ref{conj:canonical direct summand}).

The counterexample in Theorem \ref{thm:counterexample m=1} is inspired by and based on our previous work \cite{Men23}*{Example 4.3}. Let $C$ be an elliptic curve and $E$ a rank $2$ vector bundle on $C$ which is the unique non-split extension
\[
0\longrightarrow\mathcal{O}_C\longrightarrow E\longrightarrow\mathcal{O}_C\longrightarrow0.
\]
In \cite{Men23}*{Example 4.3}, we have already shown that for any integer $m\ge1$, the $m$-th symmetric power $\Sym^mE$ of $E$ does not admit any Chen--Jiang decomposition and is not globally generated after pullback by any isogeny. This is the starting point of the counterexamples in Theorems \ref{thm:counterexample m=1} and \ref{thm:counterexample m>=1}. Conjecture \ref{conj:global generation} holds for the counterexample in Theorem \ref{thm:counterexample m>=1}, which again implies that these conjectures are not equivalent in general.

\begin{ac}
{I would like to thank Mihnea Popa for bringing me into this field in $2019$, introducing me ideas and techniques around this topic, and discussing and commenting the mathematics since then. I would like to thank James M\textsuperscript{c}Kernan for helpful discussions and comments. I am partially supported by the Simons Collaboration Grant on Moduli of Varieties. I would like to thank an AI tool for language editing and reference searching. After I had discovered the counterexamples in Theorems \ref{thm:counterexample m=1} and \ref{thm:counterexample m>=1}, I learned on arXiv that similar counterexamples are also independently discovered by Houari Benammar Ammar recently in \cite{Amm26}. I would like to thank him for effective and friendly communication.}
\end{ac}

\section{Preliminaries}\label{sec:preliminaries}

We work over the field of complex numbers $\C$. Unless otherwise specified, all pairs are projective throughout this paper. A \emph{fibration} is a projective surjective morphism with connected fibers. A \emph{birational contraction} is a birational map whose inverse does not contract any divisor. For the definitions and basic results on the singularities of pairs, we refer to \cite{KM98}. We always ask the boundary $\Delta$ in a pair $(X,\Delta)$ to be effective.

We recall several definitions first.

\begin{defi}\label{defi:canonical direct summand}
A coherent sheaf $\mathcal{F}$ on a variety $Z$ is a \emph{canonical direct summand} if it is isomorphic to a direct summand of $h_*\omega_V$, where $h:V\to Z$ is projective and $V$ is smooth and possibly disconnected.
\end{defi}

\begin{defi}
Let $A$ be an abelian variety. A \emph{torsion subvariety} of $A$ is a translate of an abelian subvariety of $A$ by a torsion point which is a closed point of finite order in $A$.
\end{defi}

\begin{defi}\label{csl}
Let $\mathcal{F}$ be a coherent sheaf on an abelian variety $A$. The \emph{cohomological support loci} $V_l^i(A, \mathcal{F})$ for $i\in\N$ and $l\in\N$ are defined by
\[
V_l^i(A, \mathcal{F})=\{\alpha\in\Pic^0(A)\mid\dim H^i(A, \mathcal{F}\otimes\alpha)\geq l\}.
\]
We use $V^i(A, \mathcal{F})$ to denote $V_1^i(A, \mathcal{F})$.
\end{defi}

\begin{defi}\label{GVM}
A coherent sheaf $\mathcal{F}$ on an abelian variety $A$
\begin{enumerate}
	\item[(1)] is a GV-\emph{sheaf} if $\codim_{\Pic^0 (A)} V^i (A, \mathcal{F}) \ge i$ for every $i>0$.
	\item[(2)] is \emph{M-regular} if $\codim_{\Pic^0 (A)} V^i (A, \mathcal{F}) > i$ for every $i>0$.
	\item[(3)] \emph{satisfies} $\IT$ if $V^i (A, \mathcal{F})=\emptyset$ for every $i>0$.
\end{enumerate}
\end{defi}

It is known that M-regular sheaves are ample by \cite{Deb06}*{Corollary 3.2}, and GV-sheaves are nef by \cite{PP11b}*{Theorem 4.1}. A line bundle on $A$ is a $\GV$-sheaf if and only if it is nef by \cite{Hac04}*{Theorem 1.2} and \cite{PP11b}*{Theorem 4.1}. A line bundle on $A$ satisfies $\IT$ if and only if it is M-regular if and only if it is ample by \cite{Deb06}*{Corollary 3.2}.

\begin{defi}\label{defi:CJ decomposition}
Let $\mathcal{F}$ be a coherent sheaf on an abelian variety $A$. The sheaf $\mathcal{F}$ is said to have the \emph{Chen--Jiang decomposition} if $\mathcal{F}$ admits a finite direct sum decomposition
\[
\mathcal{F}\cong \bigoplus_{i\in I}(\alpha_i\otimes p_i^*\mathcal{F}_i),
\]
where each $A_i$ is an abelian variety, each $p_i\colon A\to A_i$ is a fibration, each $\mathcal{F}_i$ is a nonzero M-regular coherent sheaf on $A_i$, and each $\alpha_i\in\Pic^0(A)$ is a torsion line bundle.
\end{defi}

\begin{defi}\label{defi:continuous evaluation morphism}
Let $\mathcal{F}$ be a coherent sheaf on an abelian variety $A$. The \emph{continuous evaluation morphism} associated to the coherent sheaf $\mathcal{F}$ is defined by
\[
e_{\mathcal{F}}\colon \bigoplus_{\alpha\in \Pic^0(A)_{\tors}}H^0(A, \mathcal{F}\otimes\alpha)\otimes\alpha^{-1}\to \mathcal{F}
\]
which is induced from the evaluation morphisms. We define $\mathsf T(\mathcal F):=\Image e_{\mathcal{F}}$, then $\mathsf T(\mathcal F)$ is globally generated after pullback by an isogeny since $\cF$ is coherent.
\end{defi}

For an isogeny $\varphi:A'\to A$, it is not hard to see that
\begin{equation}\label{eq:isogeny}
V^i(A',\varphi^*\mathcal F)=\hat\varphi\bigl(V^i(A,\mathcal F)\bigr),\quad \mathrm{and} \quad \mathsf T(\varphi^*\mathcal F)=\varphi^*\mathsf T(\mathcal F),
\end{equation}
where $\hat\varphi$ is the dual isogeny of $\varphi$. By \cite{Men21}*{Lemma 2.9}, $\mathcal F$ is globally generated after pullback by an isogeny if and only if $\mathsf T(\mathcal F)=\mathcal F$. Since $\mathsf T(\mathcal F)$ is globally generated after pullback by an isogeny, we have that $\mathsf T(\mathsf T(\mathcal F))=\mathsf T(\mathcal F)$. The definition also gives that $H^0(A,\mathsf T(\mathcal F)\otimes\alpha)=H^0(A,\mathcal F\otimes\alpha)$ for any $\alpha\in \Pic^0(A)_{\tors}$. If $\cF$ is a canonical direct summand on an abelian variety $A$, then $\mathsf T(\mathcal F)=\mathcal F$ by \cite{PPS17}*{Theorem A}, \cite{LPS20}*{Proposition 3.6} and the proof of \cite{LPS20}*{Theorem 5.1}.

\section{Positivity of pushforwards of lc pairs}\label{sec:Positivity of pushforwards of lc pairs}

In this section, we prove Theorem \ref{thm:gg for lc}. We start with a simple lemma.

\begin{lemma}\label{lem:modification}
Let $(X,\Delta)$ be an lc pair, $m\ge1$ a rational number and $D$ a Cartier divisor on $X$ such that $D\sim_{\Q}m(K_X+\Delta)$. For a log resolution $\mu:Y\to X$ of $(X,\Delta)$, write $K_Y+\Theta=\mu^*(K_X+\Delta)$ and let
\[
J=\lceil m\Theta^-\rceil,\qquad D_Y=\mu^*D+J,\qquad \Delta_Y=\Theta+J/m.
\]
Then $(Y,\Delta_Y)$ is log smooth lc, $D_Y\sim_{\Q}m(K_Y+\Delta_Y)$, and $\mu_*\cO_Y(D_Y)=\cO_X(D)$. Here $\Theta^-$ denotes the effective negative part of $\Theta$.
\end{lemma}

\begin{proof}
By definition, $J$ is effective and exceptional, thus $\mu_*\cO_Y(D_Y)=\cO_X(D)$. By direct calculations, $(Y,\Delta_Y)$ is log smooth lc and $D_Y\sim_{\Q}m(K_Y+\Delta_Y)$.
\end{proof}

Next, we prove a result on extension from the lc locus in our relevant setting using a result on extension from the non-lc locus in \cite{Amb14}.

\begin{prop}\label{prop:lifting from the lc locus}
Let $f$ be a morphism from an lc pair $(X, \Delta)$ to an abelian variety $A$, $m>1$ a rational number and $D$ a Cartier divisor on $X$ such that $D\sim_{\Q}m(K_X+\Delta)$. Let $Z$ be a reduced union of lc centers, $\cK=f_*(\cI_Z\otimes\cO_X(D))$ and $\cQ=f_*\cO_X(D)/\cK$. If $\mathsf T(\cK)=\cK$, then the natural map
\[
 H^0(A,f_*\cO_X(D)\otimes\alpha)\longrightarrow H^0(A,\cQ\otimes\alpha)
\]
is surjective for every torsion line bundle $\alpha\in\Pic^0(A)$.
\end{prop}

\begin{proof}
If $\cK=0$, then $\cQ=f_*\cO_X(D)$, and the assertion is trivial. Thus we assume that $\cK\ne0$. Similarly, we assume that $\cQ\ne0$. Thus $Z$ is not empty. By $\mathsf T(\cK)=\cK$, we can choose finitely many torsion line bundles $\alpha_i\in\Pic^0(A)$ and sections
\[
 s_i\in H^0(A,\cK\otimes\alpha_i)=H^0(X,\cI_Z\otimes\cO_X(D)\otimes f^*\alpha_i)
\]
whose evaluation maps jointly generate $\cK$. Define the ideal sheaf $\mathfrak b$ on $X$ by
\[
 \mathfrak b\otimes\cO_X(D)=\Image\bigl(f^*\cK\longrightarrow\cO_X(D)\bigr).
\]
The sections $s_i$ generate this image after pullback by an isogeny, and $\mathfrak b\subseteq\cI_Z$.

For every irreducible component $Z_j$ of $Z$, choose an lc place $E_j$ with lc center $Z_j$. Take a log resolution $\mu:Y\to X$ which contains these divisors and principalizes $\mathfrak b$, and write
\[
 K_Y+\Theta=\mu^*(K_X+\Delta),\qquad\mathfrak b\cO_Y=\cO_Y(-F).
\]
Here $F$ is an effective divisor, every coefficient of $\Theta$ is at most one, and the support of $\Theta$, $F$ and the exceptional divisors has simple normal crossings. Choose a positive integer $n$ divisible by the orders of all $\alpha_i$, and trivialize their $n$-th powers. The sections $s_i^n$ then belong to $H^0(X,\cO_X(nD))$, and their pullbacks after $\mu$, after removing $nF$, globally generate $\cO_Y(n\mu^*D-nF)$. Take a general member of the base point free linear system $H\in|n\mu^*D-nF|$ such that $H$ is smooth, has no common components with $\Theta$, $F$ or the exceptional divisors, and meets them transversely. Let $C\in|nD|$ such that $\mu^*C=nF+H$. Let
\[
 R:=C/n\sim_{\Q}D,\qquad c:=\frac{m-1}{m}.
\]
Then $R\ge0$, $0<c<1$, and $\mu^*R=F+H/n$.

Next, we use the scheme structure on the non-lc locus defined in \cite{Amb14}. Let $\Gamma:=\Theta+cF+(c/n)H=\sum_E\gamma_E E$, then $K_Y+\Gamma=\mu^*(K_X+\Delta+cR)$. Let
\[
N(\Gamma):=\sum_{\gamma_E<0}\lfloor\gamma_E\rfloor E+\sum_{\gamma_E>1}(\lceil\gamma_E\rceil-1)E.
\]
By \cite{Amb14}*{Definition 4.3}, the ideal sheaf of the non-lc subscheme of $(X,\Delta+cR)$ is
\[
 \cI:=\mu_*\cO_Y(-N(\Gamma))\subseteq\cO_X.
\]
In particular, the part of $\Gamma$ with coefficients in $[0,1]$ contribute nothing, while the negative part contribute. We claim that $N(\Gamma)\le F$. If a prime divisor has multiplicity $e>0$ in $F$, its coefficient in $\Gamma$ is at most $1+ce$. When this coefficient exceeds one, its contribution to $N(\Gamma)$ is at most $\lceil1+ce\rceil-1\le e$, since $e$ is an integer and $c<1$. Otherwise its contribution is nonpositive. Along $H$, the coefficient is $c/n<1$, so the contribution is zero. Along every other prime divisor, the coefficient is at most one, and again the contribution is nonpositive. Thus
\[
 \mathfrak b\subseteq\mu_*\cO_Y(-F)\subseteq\mu_*\cO_Y(-N(\Gamma))=\cI.
\]
On the other hand, $\mathfrak b\subseteq\cI_Z$ implies that $F$ has a positive multiplicity along every chosen lc place $E_j$. Since the coefficient of $\Theta$ along $E_j$ is one and $c>0$, its coefficient in $\Gamma$ is bigger than one. Thus every $Z_j$ is contained in the non-lc locus of $(X,\Delta+cR)$. As $Z$ is reduced, this gives $\cI\subseteq\cI_Z$. So we have $\mathfrak b\subseteq\cI\subseteq\cI_Z$. Thus we have the natural sequence
\[
f^*f_*(\cI\otimes\cO_X(D))\longrightarrow f^*f_*(\cI_Z\otimes\cO_X(D))\longrightarrow\mathfrak b\otimes\cO_X(D)\longrightarrow\cI\otimes\cO_X(D),
\]
and therefore
\[
f_*(\cI\otimes\cO_X(D))\longrightarrow f_*(\cI_Z\otimes\cO_X(D))\longrightarrow f_*(\mathfrak b\otimes\cO_X(D))\longrightarrow f_*(\cI\otimes\cO_X(D)).
\]
Thus
\[
 f_*(\cI\otimes\cO_X(D))=f_*(\cI_Z\otimes\cO_X(D))=\cK.
\]
Pushing forward the restriction sequence gives a natural injection
\[
 \cQ\hookrightarrow f_*\bigl((\cO_X/\cI)\otimes\cO_X(D)\bigr).
\]
Finally, fix any torsion line bundle $\alpha\in\Pic^0(A)$. Since
\[
 K_X+\Delta+cR\sim_{\Q}\frac1mD+\frac{m-1}{m}D=D
\]
and $f^*\alpha$ is $\Q$-linearly trivial, by applying \cite{Amb14}*{Theorem 6.2} to the line bundle $\cO_X(D)\otimes f^*\alpha$ we get a surjection
\[
 H^0(X,\cO_X(D)\otimes f^*\alpha)\longrightarrow H^0\bigl(X,(\cO_X/\cI)\otimes\cO_X(D)\otimes f^*\alpha\bigr).
\]
Since the above map factors through $H^0(A,\cQ\otimes\alpha)$ which injects into the target of the above map, we conclude that
\[
H^0(X,\cO_X(D)\otimes f^*\alpha)\longrightarrow H^0(A,\cQ\otimes\alpha)
\]
is surjective, and we are done.
\end{proof}

Next, we study the positivity of the restricted image of a particular sequence of coherent sheaves, which appears naturally for our purposes. We first study the case when $m=1$. The situation is much simpler if the boundary divisor of the pair has no components with coefficient one, but requires essential work otherwise. In particular, input from Viehweg's cyclic covering trick and the theory of mixed Hodge modules is needed here.

\begin{prop}\label{prop:gg for lc restricted image m=1}
Let $(Y,S+B)$ be a log smooth lc pair with $S=\lfloor S+B\rfloor$, $h:Y\to A$ a morphism to an abelian variety $A$, and $T$ a reduced union of minimal lc centers of the same dimension. If $D\sim_{\Q}K_Y+S+B$ is a Cartier divisor on $Y$, then
\[
 \cR:=\Image\bigl(h_*\cO_Y(D)\longrightarrow (h|_T)_*\cO_T(D|_T)\bigr)
\]
is a canonical direct summand.
\end{prop}

\begin{proof}
Let $M:=D-K_Y-S$. Since $M\sim_{\Q}B$, choose an integer $N\ge2$ such that $NB$ is integral and $NM\sim NB$, and let $s\in H^0(Y,\cO_Y(NM))$ be a section with zero divisor $NB$. Let $\pi:Z\to Y$ be the normalization of the cyclic cover $\Spec_Y\bigl(\oplus_{i=0}^{N-1}\cO_Y(-iM)\bigr)$, where multiplication is defined by $s$. The cover may be disconnected. By \cite{EV92}*{Corollary 3.11}, $\pi$ is finite and flat, and
\[
 \pi_*\cO_Z=\bigoplus_{i=0}^{N-1}\cO_Y(-iM+\lfloor iB\rfloor).
\]
Moreover, $Z$ has quotient singularities, hence rational singularities by \cite{EV92}*{Lemma 3.24 and the proof of Claim 5.13}. Duality for a finite flat morphism gives the corresponding decomposition
\[
 \pi_*\omega_Z=\bigoplus_{i=0}^{N-1}\omega_Y\otimes\cO_Y(iM-\lfloor iB\rfloor).
\]

Take a resolution $\rho:X\to Z$ functorial for smooth morphisms as in \cite{BM08}*{Theorem 1.1}, and let
\[
 q:=\pi\rho,\qquad E:=q^*S,\qquad T':=q^{-1}(T).
\]
\'Etale locally around each point of $Y$, choose coordinates in which the components of $\Supp(S+B)$ are coordinate hyperplanes. By further shrinking the \'etale neighborhood, we can extract an $N$-th root of the local unit in the equation of the cyclic covering and make the equation only involve the coordinates defining $B$. The coordinates defining $S$ do not occur in the equation. Thus the cyclic cover is \'etale locally a product with one factor the coordinate space given by the coordinates of $S$. Normalization preserves this product structure. Since $\rho$ is functorial for smooth morphisms, this product structure is preserved, and the natural $\mu_N$-action on the cyclic cover extends to its normalization and lifts to the functorial resolution. Thus $E$ is reduced with simple normal crossing support, and $T'$ is a reduced union of minimal lc centers of $(X,E)$. Note that minimal lc centers of a log smooth lc pair are smooth and pairwise disjoint. The scheme $\pi^{-1}(T)$ with the reduced scheme structure is the normalized cyclic cover of $T$ defined by $s|_T$, and $\rho|_{T'}$ is a resolution of it. Moreover, $\pi^{-1}(T)$ also has quotient and rational singularities. Thus $\rho_*\omega_X=\omega_Z$, and the analogous equality holds for $\rho|_{T'}$. Thus we have
\[
 q_*\omega_X(E)=\bigoplus_{i=0}^{N-1}\omega_Y(S)\otimes\cO_Y(iM-\lfloor iB\rfloor),\quad(q|_{T'})_*\omega_{T'}=\bigoplus_{i=0}^{N-1}\omega_T\otimes\cO_T(iM|_T-\lfloor iB\rfloor|_T).
\]
No components of $B$ contain a component of $T$, and the simple normal crossing condition implies $\lfloor iB\rfloor|_T=\lfloor iB|_T\rfloor$. The minimality condition on lc centers gives $\omega_Y(S)|_T\cong\omega_T$.

Fix an ordering of the components of $S$, and use the induced ordering to take iterated residues along the components of $E$ which define each component of $T'$. These residue maps commute with the natural $\mu_N$-action and preserve the indexed direct summands above. Over $T\setminus\Supp B$, $\pi$ is \'etale, so the map on each direct summand is the ordinary restriction map given by adjunction. Thus the equality holds over the entire $T$, since the target of each direct summand is a line bundle on $T$. By $\lfloor B\rfloor=0$ and $\omega_Y(S)|_T\cong\omega_T$, the residue map at $i=1$ is exactly the natural restriction map
\[
 \cO_Y(D)\longrightarrow\cO_T(D|_T).
\]
After applying $h_*$, it follows that $\cR$ is the $i=1$ direct summand of
\[
\cR':=\Image\bigl((hq)_*\omega_X(E)\longrightarrow ((hq)|_{T'})_*\omega_{T'}\bigr).
\]

Next, we show that $\cR'$ is a direct summand of $((hq)|_{T'})_*\omega_{T'}$. Let $j:X\setminus E\rightarrow X$ and $i:T'\rightarrow X$ be the inclusions. By iterating the residue sequence \cite{FFS14}*{(2.4.3)} for the irreducible components defining $T'$, we get a morphism of mixed Hodge modules
\[
 j_*\mathbb D\mathbb Q^H_{X\setminus E}[-\dim X]\longrightarrow i_*\mathbb D\mathbb Q^H_{T'}[-\dim T'],
\]
where $\mathbb Q^H$ denotes the constant Hodge complex and $\mathbb D$ denotes the duality functor. For the underlying filtered right $\mathcal D$-modules, their lowest nonzero graded pieces occur at index zero and are $\omega_X(E)$ and $i_*\omega_{T'}$ respectively, and the morphism between them is exactly the iterated residue map
\[
\omega_X(E)\longrightarrow i_*\omega_{T'}
\]
by \cite{FFS14}*{Theorems 1 and 2.3}. By applying the zeroth cohomological direct image under $hq$ to the above morphism of mixed Hodge modules and using \cite{FFS14}*{Theorem 1 and (1.1.3)}, we get a morphism of mixed Hodge modules 
\[
\cM\longrightarrow\cN
\]
whose zeroth graded piece is exactly the map which defines $\cR'$. Since every morphism of mixed Hodge modules is strict with respect to the Hodge filtration, $\cR'$ is the zeroth graded piece of the image of $\cM\longrightarrow\cN$. Since $T'$ is equidimensional and smooth, $\cN$ is polarizable pure of weight $-\dim T'$ by the direct image theorem \cite{Sai17}*{Theorem 1.4}. Thus semisimplicity implies that the image of $\cM\longrightarrow\cN$ splits and thus its zeroth graded piece also splits. Thus we have a projector
\[
P:((hq)|_{T'})_*\omega_{T'}\longrightarrow ((hq)|_{T'})_*\omega_{T'}
\]
given by the splitting such that $P^2=P$ and $\Image P=\cR'$. Since $\cR$ is a direct summand of $\cR'$ and $\cR'$ is a direct summand of $((hq)|_{T'})_*\omega_{T'}$, the canonical direct summand assertion follows.

Moreover, we claim that we can choose a $\mu_N$-equivariant projector from $((hq)|_{T'})_*\omega_{T'}$ onto $\cR$. We will need this later. By construction, the morphism $\cM\to\cN$ is $\mu_N$-equivariant. Thus its induced morphism on the zeroth graded pieces is $\mu_N$-equivariant, so its image $\cR'$ is $\mu_N$-stable. Thus $\cR$ is also $\mu_N$-stable by construction. Choose a projector from $((hq)|_{T'})_*\omega_{T'}$ onto $\cR$ whose existence has already been established. Every conjugate of this projector under $\mu_N$ maps onto $\cR$ and restricts to the identity on $\cR$. By taking an average of these conjugates of the projector, we get a $\mu_N$-equivariant projector from $((hq)|_{T'})_*\omega_{T'}$ onto $\cR$.
\end{proof}

Next, we study the similar problem as in Proposition \ref{prop:gg for lc restricted image m=1} in the case when $m>1$.

\begin{prop}\label{prop:gg for lc restricted image m>1}
Let $(X,S+B)$ be a log smooth lc pair with $S=\lfloor S+B\rfloor$, $f:X\to A$ a morphism to an abelian variety $A$, $m>1$ a rational number and $D$ a Cartier divisor on $X$ such that $D\sim_{\Q}m(K_X+S+B)$. Let $T$ be a reduced union of minimal lc centers of the same dimension and
\[
 \mathcal Q:=\Image\bigl(f_*\mathcal O_X(D)\longrightarrow(f|_T)_*\mathcal O_T(D|_T)\bigr).
\]
Assume that the natural map
\[
 H^0(A,f_*\cO_X(D)\otimes\alpha)\longrightarrow H^0(A,\cQ\otimes\alpha)
\]
is surjective for every torsion line bundle $\alpha\in\Pic^0(A)$. Then $\mathsf T(\cQ)$ is a canonical direct summand. 
\end{prop}

\begin{proof}
We can assume that $\mathsf T(\cQ)\ne0$. Choose torsion line bundles $\alpha_1,\ldots,\alpha_r\in\Pic^0(A)$ and sections $q_i\in H^0(A,\cQ\otimes\alpha_i)$ whose evaluation maps jointly generate $\mathsf T(\cQ)$. By the lifting assumption and the projection formula, we can choose the lifted sections $u_i\in H^0(X,\cO_X(D)\otimes f^*\alpha_i)$ of each $q_i$. Let $\mathfrak b\subseteq\cO_X$ be their common base ideal defined by
\[
 \mathfrak b\otimes\cO_X(D)=\Image\bigl(\bigoplus_{i=1}^r f^*\alpha_i^{-1}\xrightarrow{(u_i)}\cO_X(D)\bigr).
\]
This ideal is nonzero because $\mathsf T(\cQ)\ne0$.

Call an irreducible component $T_j$ of $T$ \emph{retained} if $u_i|_{T_j}$ is nonzero for some $i$. Thus $\mathfrak b=\cO_X$ at the generic point of every retained center. For each $T_j$ which is not retained, first blow up $T_j$ if it has codimension at least two, and denote the exceptional divisor by $E_j$. If $T_j$ is a divisor, take $E_j=T_j$. The minimal lc centers are smooth and pairwise disjoint, so these operations do not affect the retained centers. Each divisor $E_j$ has coefficient one in the crepant pullback of $S+B$, and $\ord_{E_j}(\mathfrak b)>0$. 

Next, principalize the pullback of $\mathfrak b$ by blowups with smooth centers which are in its cosupport and have normal crossings with the boundary and the exceptional divisors, as in \cite{BM08}*{Theorem 1.3}. Denote the composition of them by $\mu:Y\to X$, and write
\[
 K_Y+\Theta=\mu^*(K_X+S+B),\qquad \mathfrak b\cO_Y=\cO_Y(-F_0).
\]
By abuse of notation, we still use the notation $E_j$ for the strict transforms of the divisors extracted from $T_j$ above. We can also assume that the union of $F_0$, strict transform of $S+B$ and the exceptional divisors has simple normal crossing support. Apply Lemma \ref{lem:modification}, and let
\[
 J:=\lceil m\Theta^-\rceil,\qquad D_Y:=\mu^*D+J,\qquad \Delta_Y:=\Theta+\frac{J}{m}.
\]
Then $(Y,\Delta_Y)$ is log smooth lc,
\[
 D_Y\sim_{\Q}m(K_Y+\Delta_Y),\qquad \mu_*\cO_Y(D_Y)=\cO_X(D).
\]
Let $F:=F_0+J$. After removing $F_0$, the pullbacks of the $u_i$ by $\mu$ generate $\cO_Y(D_Y-F)$ with their respective torsion twists.

We claim that the strict transforms of the retained centers are still minimal lc centers of $(Y,\Delta_Y)$. Since $\mathfrak b=\cO_X$ at the generic point of every retained center, the blowups we choose are isomorphisms at these generic points. If the strict transform of a retained center $T_j$ is not a minimal lc center of $(Y,\Delta_Y)$, then it is not a minimal lc center of $(Y,\Theta)$ since the term $J/m$ changes the negative coefficients of $\Theta$ into $[0,1/m)$ and thus it does not create new divisors with coefficient one. Then at some intermediate blowup, we have created the very first exceptional divisor with coefficient one in the crepant boundary such that it intersects the strict transform of $T_j$. Blowing up a smooth center $C$ having normal crossings with the crepant boundary $\Gamma$ gives the coefficient of the exceptional divisor
\[
 \sum_{\Gamma_k\supset C}\coef_{\Gamma_k}(\Gamma)-\codim C+1.
\]
in the new crepant boundary. Since all crepant boundary coefficients are at most one, this equals one only when $C$ is an lc center. Any lc center meeting a retained minimal lc center contains that retained minimal lc center. Since $T_j$ is not in the cosupport of $\mathfrak b$, we deduce that $C$ is not in the cosupport of the transform of $\mathfrak b$, which contradicts our previous choice of $C$. Thus our claim is proved. Thus the reduced union $T_Y$ of the strict transforms of the retained minimal lc centers consists of minimal lc centers of $(Y,\Delta_Y)$, all of the original dimension. Moreover, $\mu$ is an isomorphism at their generic points, and $F$ does not contain any irreducible component of $T_Y$.

Choose an integer $n\ge2$ divisible by the orders of all the $\alpha_i$. After removing $nF_0$, the $n$-th powers of the pullbacks of the $u_i$ by $\mu$ generate $\cO_Y(n(D_Y-F))$. Choose a general member $H$ in $|n(D_Y-F)|$. By Bertini's theorem, the support of $H$, $\Delta_Y+F$, and the exceptional divisors has simple normal crossings, and $H$ has no common irreducible components with $\Delta_Y+F$ or exceptional divisors. It also contains no irreducible components of $T_Y$. Let
\[
 R:=F+\frac{H}{n}\sim_{\Q}D_Y, \qquad c:=\frac{m-1}{m}.
\]
Let $S_0$ be the sum of the coefficient one irreducible components of $\Delta_Y$ which are not in $F$, and define
\[
 V:=\lfloor\Delta_Y+cR\rfloor-S_0=\lfloor\Delta_Y+cF\rfloor-S_0.
\]
The second equality above follows from the choice of $H$ and $c/n<1$. We have $0\le V\le F$. Indeed, if a prime divisor has coefficient $e\ge1$ in $F$ and coefficient $b\in[0,1]$ in $\Delta_Y$, then $0\le\lfloor b+ce\rfloor\le e$. If $e=0$, it is not in $V$. We also have
\[
 D_Y-V\sim_{\Q}K_Y+S_0+\{\Delta_Y+cR\}.
\]
This is a log smooth lc boundary whose reduced part is $S_0$. Since $F$ does not contain any irreducible component of $T_Y$, every coefficient one irreducible component of $\Delta_Y$ defining an irreducible component of $T_Y$ is not in $F$ and thus in $S_0$. Therefore the irreducible components of $T_Y$ are still minimal lc centers of this new pair. By Proposition \ref{prop:gg for lc restricted image m=1}, the sheaf
\[
 \cR=\Image\bigl((f\mu)_*\cO_Y(D_Y-V)\longrightarrow((f\mu)|_{T_Y})_*\cO_{T_Y}((D_Y-V)|_{T_Y})\bigr)
\]
is a canonical direct summand.

Next, we identify $\cR$ with a subsheaf of $\cQ$. Pushing forward the inclusion $\cO_Y(D_Y-V)\rightarrow\cO_Y(D_Y)$ by $f\mu$ gives another map with source $(f\mu)_*\cO_Y(D_Y-V)$
\[
 (f\mu)_*\cO_Y(D_Y-V)\longrightarrow f_*\cO_X(D)\longrightarrow\cQ.
\]
We claim that this map and the restriction map defining $\cR$ have the same kernel. Consider a local section of their source and its descended section of $\cO_X(D)$. Over each retained center $T_j$, $\mu$ is generically an isomorphism and neither $J$ nor $V$ contains its strict transform (note that $0\le V\le F$). Thus the local section vanishes on the strict transform of $T_j$ if and only if its descended section vanishes on $T_j$. Note that this can be checked at the generic points, since the sections which we are considering here are sections of line bundles on the disjoint union of integral minimal lc centers.

For a center $T_j$ which is not retained, the divisor $E_j$ has coefficient one in $\Theta$, coefficient zero in $J$, and positive coefficient in $F$. It therefore has positive coefficient in $V$. Every section of $\cO_Y(D_Y-V)$ descends to a section whose pullback vanishes along $E_j$, since $E_j$ has coefficient zero in $J$. Such a descended section vanishes on $T_j$. Thus the descended section always vanishes on the centers which are not retained. Thus we have proved the equality of kernels, and it gives a natural injection $\cR\hookrightarrow\cQ$. We have the following commutative diagram.
\begin{center}
\begin{tikzcd}
			(f\mu)_*\cO_Y(D_Y-V) \arrow[r, hook,""] \arrow[d, two heads,""] & f_*\cO_X(D) \arrow[d, two heads,""] \\
			 \cR \arrow[r, hook,""] & \cQ 
\end{tikzcd}
\end{center}

Finally, $V\le F$ implies that every section $\mu^*u_i$ is in $H^0(Y,\cO_Y(D_Y-V)\otimes(f\mu)^*\alpha_i)$. Under the injection $\cR\hookrightarrow\cQ$, its restriction to $T_Y$ maps to $q_i$. Since the evaluation maps of the $q_i$ jointly generate $\mathsf T(\cQ)$, we obtain
\[
 \mathsf T(\cQ)\subseteq\cR.
\]
Conversely, since $\cR$ is a canonical direct summand on an abelian variety $A$, we have $\mathsf T(\cR)=\cR$. Since $\cR\hookrightarrow\cQ$ is injective, we have
\[
 \cR=\mathsf T(\cR)\subseteq\mathsf T(\cQ).
\]
Thus $\mathsf T(\cQ)=\cR$, and $\mathsf T(\cQ)$ is a canonical direct summand on $A$.
\end{proof}

The following theorem is central to the proof of our main theorem. The theory of singular Hermitian metrics with semipositive curvature and the minimal extension property plays a crucial role here, see \cite{HPS18} for related definitions and results regarding the theory.

\begin{thm}\label{thm:image splitting}
Let $(X,S+B)$ be a log smooth lc pair with $S=\lfloor S+B\rfloor$, $f:X\to A$ a morphism to an abelian variety $A$, $m\ge1$ a rational number and $D$ a Cartier divisor on $X$ such that $D\sim_{\Q}m(K_X+S+B)$. Let $T$ be a reduced union of minimal lc centers of the same dimension such that each of them dominates $A$. Then the sheaf
\[
\mathcal Q:=\Image\bigl(f_*\mathcal O_X(D)\longrightarrow(f|_T)_*\mathcal O_T(D|_T)\bigr)
\]
admits a singular Hermitian metric with semipositive curvature and the minimal extension property. In particular, every nonzero morphism $\cQ\to\cO_A$ splits.
\end{thm}

\begin{proof}
We may assume that $\cQ\ne0$ and $\dim A>0$. Note that $f$ is surjective and $\cQ$ is torsion-free by the assumptions. If $m=1$, $\cQ$ is a canonical direct summand by Proposition \ref{prop:gg for lc restricted image m=1}. Since $\cQ$ is torsion-free, $\cQ$ is isomorphic to a direct summand of the pushforward of the canonical bundle under a projective morphism which dominates $A$ when restricted to each irreducible component of the source. By \cite{HPS18}*{Theorem 21.1}, the pushforward of the canonical bundle of each irreducible component of the source admits a singular Hermitian metric with semipositive curvature and the minimal extension property. Then the direct sum of all of them with the induced metric also admits a singular Hermitian metric with semipositive curvature and the minimal extension property by definition. Moreover, consider a torsion-free quotient sheaf of a torsion-free coherent sheaf admitting a singular Hermitian metric with semipositive curvature and the minimal extension property. By definition, the induced quotient metric on the quotient sheaf satisfies the same properties, see e.g. \cite{LS24}*{Proposition 2.3}. In summary, $\cQ$ admits a singular Hermitian metric with semipositive curvature and the minimal extension property.

Next, we assume that $m>1$. Let $\cF:=f_*\cO_X(D)$. We apply the principalization argument in the proof of Proposition \ref{prop:gg for lc restricted image m>1} to the relative base ideal defined by
\[
 \mathfrak b\otimes\cO_X(D)=\Image(f^*\cF\longrightarrow\cO_X(D)).
\]
Discard the irreducible components of $T$ on which every local section of $\cF$ restricts to zero. This leaves $\cQ$ unchanged. At the generic point of each \emph{retained} irreducible component, $\mathfrak b=\cO_X$. Choose the principalization $\mu:Y\to X$ of $\mathfrak b$ as in that proof and write
\[
 \begin{gathered}
 K_Y+\Theta=\mu^*(K_X+S+B),\qquad\mathfrak b\cO_Y=\cO_Y(-F_0),\qquad J:=\lceil m\Theta^-\rceil,\\
 D_Y:=\mu^*D+J,\qquad \Delta_Y:=\Theta+J/m,\qquad h:=f\mu,\qquad F:=F_0+J,\qquad P:=D_Y-F.
 \end{gathered}
\]
We have $h_*\cO_Y(D_Y)=\cF$ and $h^*\cF\twoheadrightarrow\cO_Y(P)$. By the same argument, the reduced union $T_Y$ of the strict transforms of the retained lc centers consists of minimal lc centers of the same dimension of $(Y,\Delta_Y)$, and no irreducible component of $T_Y$ is contained in $F$. Let $S_0$ be the sum of the coefficient one irreducible components of $\Delta_Y$ which are not in $F$, and define
\[
 c:=\frac{m-1}{m},\qquad V:=\lfloor\Delta_Y+cF\rfloor-S_0,\qquad B_1:=\{\Delta_Y+cF\},\qquad D':=D_Y-V.
\]
For $0\le b\le1$ and an integer $e\ge1$, we have that $0\le\lfloor b+ce\rfloor\le e$ and thus $0\le V\le F$. Since $h_*\cO_Y(P)=\cF$ and $0\le V\le F$, we deduce that $h_*\cO_Y(D')=\cF$. Moreover,
\[
 D'\sim_{\Q}K_Y+S_0+B_1+cP.
\]
Let $g:=h|_{T_Y}$, $D'_T:=D'|_{T_Y}$, $P_T:=P|_{T_Y}$, and $B_T:=B_1|_{T_Y}$. By adjunction, we have that $D'_T\sim_{\Q}K_{T_Y}+B_T+cP_T$ and $(T_Y,B_T)$ is log smooth klt. Since $\mu$ is generically an isomorphism along each retained center and neither $J$ nor $V$ contains its strict transform, the maps $h_*\cO_Y(D')\to g_*\cO_{T_Y}(D'_T)$ and $h_*\cO_Y(D')\to f_*\cO_X(D)\to\cQ$ have the same kernel. Thus we identify
\[
 \cQ\subseteq\cE:=g_*\cO_{T_Y}(D'_T).
\]

\emph{Step 1.} Recall that $\cF=h_*\cO_Y(P)=h_*\mathcal O_Y(D')=h_*\mathcal O_Y(D_Y)$. Denote $\cK:=\ker(\mathcal F\twoheadrightarrow\mathcal Q)$. We now prove that $g^*\mathcal F\twoheadrightarrow\mathcal O_{T_Y}(P_T)$, the pullback of the relative evaluation map of $\cO_Y(P)$, factors through $g^*\mathcal Q$. The composition of this map and the injective map $\mathcal O_{T_Y}(P_T)\hookrightarrow\cO_{T_Y}(D'_T)$ given by the multiplication by the nonzero section of $(F-V)|_{T_Y}$ is $g^*\mathcal F\to\cO_{T_Y}(D'_T)$, the pullback of the relative evaluation map of $\mathcal O_Y(D')$. The multiplication map is injective since $0\leq F-V\leq F$ and no irreducible component of $T_Y$ is contained in $F$. Under the identification of $\mathcal Q$ with $\Image(\mathcal F\to g_*\mathcal O_{T_Y}(D'_T))$, the map $g^*\mathcal F\to\cO_{T_Y}(D'_T)$ vanishes on the image of $g^*\mathcal K$. Thus $g^*\mathcal F\twoheadrightarrow\mathcal O_{T_Y}(P_T)$ also vanishes on the image of $g^*\mathcal K$. By the exact sequence $g^*\mathcal K\to g^*\mathcal F\to g^*\mathcal Q\to 0$, we deduce that $g^*\mathcal F\twoheadrightarrow\mathcal O_{T_Y}(P_T)$ factors through a surjection
\[
 g^*\cQ\twoheadrightarrow\cO_{T_Y}(P_T).
\]

\emph{Step 2.} We now prove that there exists a splitting surjection $\pi:\cE\twoheadrightarrow\cQ$. Choose a sufficiently ample Cartier divisor $M$ on $A$ such that $M/c$ is Cartier and $\cF\otimes\cO_A(M/c)$ is globally generated. Since we have the natural surjection $h^*\cF\twoheadrightarrow\cO_Y(P)$, $\cO_Y(P+h^*(M/c))$ is globally generated. Choose a general member $H\in|P+h^*(M/c)|$ such that $H$ and $S_0+B_1$ have no common irreducible components and $H+S_0+B_1$ has simple normal crossing support. We have that
\[
 D'+h^*M\sim_{\Q}K_Y+S_0+B_1+cH,
\]
where $(Y,S_0+B_1+cH)$ is log smooth lc with the reduced part of its boundary $S_0$. Since no irreducible component of $T_Y$ is contained in $F$, $T_Y$ is a union of minimal lc centers of the same dimension of $(Y,S_0+B_1+cH)$. The restricted image of $D'+h^*M$ when restricted to $T_Y$ is $\cQ\otimes\cO_A(M)$ which is in $\cE\otimes\cO_A(M)$. We apply the proof of Proposition \ref{prop:gg for lc restricted image m=1} to our setting here. By restricting the equivariant projector with respect to the cyclic covering to the $i=1$ part, we get a splitting surjection $\cE\otimes\cO_A(M)\twoheadrightarrow\cQ\otimes\cO_A(M)$. By twisting it by $\cO_A(-M)$, we get the desired splitting surjection $\pi:\cE\twoheadrightarrow\cQ$.

\emph{Step 3.1.} In the following steps, we work in the analytic setting. Let $L:=\cO_{T_Y}(D'_T)\otimes\omega_{T_Y/A}^{-1}$. Fix a nonzero invariant holomorphic top form on $A$, identifying $\omega_{T_Y/A}$ with $\omega_{T_Y}$. Choose a sufficiently divisible positive integer $N$ such that $NB_T$ is integral, $Nc\in\mathbb Z$, and we fix an isomorphism
\[
 L^{\otimes N}\cong\cO_{T_Y}(NB_T)\otimes\cO_{T_Y}(P_T)^{\otimes Nc}.
\]
Let $h_{NB_T}$ be the canonical singular Hermitian metric associated with the effective divisor $NB_T$. For any metric $k$ on $\cO_{T_Y}(P_T)$, it induces the tensor product metric $h_{NB_T}\otimes k^{\otimes Nc}$ on the right-hand side. By the fixed isomorphism, $h_{NB_T}\otimes k^{\otimes Nc}$ induces a metric on the left-hand side. By taking its $N$-th root, we get a metric on $L$ denoted by $h_{B_T}k^c$. The similar construction works when $k$ is only defined over an open subset of $T_Y$. By abuse of notation, we still use $h_{B_T}k^c$ to denote the outcome of the construction.

Let $T_{Y,a}:=g^{-1}(a)$. Choose an open subset $A^\circ\subseteq A$ over which $g$, $h$, all the strata of $\Supp B_T$ and all the strata of $\Supp(\Delta_Y+F)$ are smooth, $h^0(T_{Y,a},\cO_{T_Y}(D'_T)|_{T_{Y,a}})$ is constant, and $\cE$ is locally free and compatible with base change. For any open subset $U\subseteq A^\circ$ and any continuous, finite, positive definite Hermitian metric $k$ on $\cO_{T_Y}(P_T)|_{g^{-1}(U)}$, we obtain the metric $h_{B_T}k^c$ on $L|_{g^{-1}(U)}$ by the previous construction. Following \cite{HPS18}*{Section 14}, we define the fiberwise $L^2$ metric $H(k)$ on $\cE|_U$ by
\[
\|e\|_{H(k),a}^2:=\int_{T_{Y,a}}|e|_{h_{B_T}k^c}^2, \qquad a\in U,\quad e\in\cE(a)\cong H^0(T_{Y,a}, \omega_{T_{Y,a}}\otimes L|_{T_{Y,a}}),
\]
where $\cE(a):=\cE_a\otimes_{\cO_{A,a}}k(a)$ is the fiber of $\cE$ over $a$, and $|e|_{h_{B_T}k^c}^2$ is the nonnegative measurable $(\dim T_{Y,a},\dim T_{Y,a})$-form associated with $e$ and the metric $h_{B_T}k^c$. If $T_{Y,a}$ is disconnected, the integral is the sum of the integrals over its connected components. Since all the strata of $\Supp B_T$ are smooth over $A^\circ$, the coefficient of each irreducible component of $B_T$ is less than one, and the metric $k$ is assumed to be continuous, finite and positive definite, we deduce that the metric $H(k)$ is continuous, finite and positive definite by basic analysis and the dominated convergence.

By Steps 1 and 2, we have a surjection
\[
\rho:g^*\cE\xtwoheadrightarrow{g^*\pi}g^*\cQ\twoheadrightarrow\cO_{T_Y}(P_T).
\]
For $a\in U$ and $x\in T_{Y,a}$, denote its map on fibers over $x$ by
\[
 \rho_x:\cE(a)\twoheadrightarrow\cO_{T_Y}(P_T)(x).
\]
Let $\Psi(k)$ be the quotient metric induced by $g^*H(k)$ and $\rho$. Explicitly,
\[
 |v|_{\Psi(k),x}^2:=\min_{\substack{e\in\cE(a)\\ \rho_x(e)=v}}\|e\|_{H(k),a}^2,\qquad v\in\cO_{T_Y}(P_T)(x).
\]
The minimum exists since $H(k)$ induces a finite positive definite Hermitian form on the finite dimensional vector space $\cE(a)$. Since $\rho$ is surjective and $H(k)$ is continuous, finite and positive definite, the quotient metric $\Psi(k)$ is continuous, finite, and positive definite on $\cO_{T_Y}(P_T)|_{g^{-1}(U)}$. In summary, whenever $k$ is continuous, finite, and positive definite, both of the metrics $H(k)$ and $\Psi(k)$ have the same three properties. If $k$ has semipositive curvature over $g^{-1}(U)$, then so does $h_{B_T}k^c$. Note that we assume that $k$ is continuous, finite and positive definite. Then the metric $h_{B_T}k^c$ has trivial multiplier ideal on $g^{-1}(U)$, since $(T_Y,B_T)$ is klt and the weights of $k$ are locally bounded. The metric $H(k)$ on $\cE|_U$ constructed in our setting is exactly the fiberwise $L^2$ metric constructed in \cite{HPS18}*{Theorem 21.1 and Section 22}. By \cite{HPS18}*{Theorem 21.1}, the metric $H(k)$ has semipositive curvature over $U$. Thus $g^*H(k)$ and $\Psi(k)$ have semipositive curvature over $g^{-1}(U)$. When $g^{-1}(U)$ is disconnected, we apply our discussion to its connected components and take the orthogonal direct sum.

All the equalities and inequalities between metrics below compare their squared norms. The constructions of $H(k)$ and $\Psi(k)$ preserve inequalities between squared norms. For $t>0$, let $tk$ denote the metric whose squared norms are $t$ times those of $k$. Then
\[
 H(tk)=t^cH(k),\qquad\Psi(tk)=t^c\Psi(k).
\]
For each $a\in U$, the Hermitian form induced by $H(k)$ on $\cE(a)$ and the metric $\Psi(k)|_{T_{Y,a}}$ only depend on $k|_{T_{Y,a}}$. Thus for every open subset $V\subseteq U$,
\[
 H(k|_{g^{-1}(V)})=H(k)|_V,\qquad\Psi(k|_{g^{-1}(V)})=\Psi(k)|_{g^{-1}(V)}.
\]
The same formulas also define $\Psi(k)$ for any continuous, finite, positive definite metric $k$ on $\cO_{T_Y}(P_T)|_{T_{Y,a}}$, without the necessity to extend $k$ to nearby fibers.

For two continuous, finite, positive definite metrics $k$ and $l$ on $\cO_{T_Y}(P_T)|_{T_{Y,a}}$, we define
\[
 d_a(k,l):=\max_{x\in T_{Y,a}}\left|\log\frac{k}{l}(x)\right|,\qquad\frac{k}{l}(x):=\frac{|v|_{k,x}^2}{|v|_{l,x}^2},\qquad v\in\mathcal O_{T_Y}(P_T)(x)\setminus\{0\}.
\]
The ratio of the squared norms is independent of $v\in\mathcal O_{T_Y}(P_T)(x)\setminus\{0\}$. Thus we get a continuous function $\frac{k}{l}:T_{Y,a}\to(0,+\infty)$. Since $T_{Y,a}$ is compact, $d_a(k,l)$ is finite. Note that $d_a(k,l)$ is nothing but the maximum absolute difference of the weights of $k$ and $l$. Let $\delta:=d_a(k,l)$, then $e^{-\delta}l\le k\le e^\delta l$. By the properties of $\Psi$, we deduce that
\[
 e^{-c\delta}\Psi(l)\le\Psi(k)\le e^{c\delta}\Psi(l).
\]
Thus we have that
\begin{equation}\label{eq:metric iteration contraction}
 d_a\bigl(\Psi(k),\Psi(l)\bigr)\le c\,d_a(k,l).
\end{equation}

\emph{Step 3.2.} We now construct a sequence of metrics and its limit locally over $A$ by repeatedly applying $\Psi$. Fix any point of $A$ and choose a sufficiently small coordinate ball $U:=B(0,5R)$ (by abuse of notation) centered around that point over which finitely many sections $s_1,\ldots,s_{N_0}\in H^0(U,\cE)$ generate $\cE|_U$. Then the sections $\rho(g^*s_1),\ldots,\rho(g^*s_{N_0})$ generate $\cO_{T_Y}(P_T)|_{g^{-1}(U)}$. They define a smooth, finite, positive definite Hermitian metric $k_0$ with semipositive curvature on $\cO_{T_Y}(P_T)|_{g^{-1}(U)}$ as follows. If $p$ is a holomorphic frame of $\cO_{T_Y}(P_T)|_{g^{-1}(U)}$ over an open subset $W\subseteq g^{-1}(U)$ and $\rho(g^*s_i)|_{W}=f_i p$, we define
\[
 |p|_{k_0}^2:=\frac{1}{\sum_{i=1}^{N_0}|f_i|^2}.
\]
Since $\rho(g^*s_1),\ldots,\rho(g^*s_{N_0})$ generate $\cO_{T_Y}(P_T)|_{g^{-1}(U)}$, the holomorphic functions $f_1,\ldots,f_{N_0}$ have no common zeros and the denominator is positive everywhere over $W$. It is not hard to see that the definition of $k_0$ is independent of the holomorphic frame we choose and the construction glues to give a metric $k_0$ on $\cO_{T_Y}(P_T)|_{g^{-1}(U)}$. Its weight $\log(\sum_{i=1}^{N_0}|f_i|^2)$ is smooth and plurisubharmonic. Thus $k_0$ is smooth, finite, positive definite and has semipositive curvature.

Over $g^{-1}(U\cap A^\circ)$, with the initial input the restriction of $k_0$ to $g^{-1}(U\cap A^\circ)$, we define successively
\[
 k_{j+1}:=\Psi(k_j),\qquad j\in\N.
\]
By the previous discussion, every $k_j$ is continuous, finite, positive definite, and has semipositive curvature. For every compact subset $K\subseteq U\cap A^\circ$, $g^{-1}(K)$ is still compact since $g$ is proper. Thus we have the following constant which is well defined
\[
 C_K:=\max_{x\in g^{-1}(K)}\left|\log\frac{k_1}{k_0}(x)\right|<+\infty.
\]
By applying the inequality \eqref{eq:metric iteration contraction} to each fiber over $K$, we deduce that
\[
 \max_{x\in g^{-1}(K)}\left|\log\frac{k_{j+1}}{k_j}(x)\right|\le c^j C_K.
\]
Recall that $\frac{k_{j+1}}{k_j}(x)$ is independent of $v\in\mathcal O_{T_Y}(P_T)(x)\setminus\{0\}$. Since $0<c<1$, $\log\frac{k_{j+1}}{k_j}(x)$ form a uniformly Cauchy sequence over compact subsets. Thus the sum of them 
\[
\sum_{j=0}^{\infty}\log\frac{k_{j+1}}{k_j}(x)
\]
exists as a finite and continuous function and is independent of $v\in\mathcal O_{T_Y}(P_T)(x)\setminus\{0\}$. Thus $k_0 e^{\sum_{j=0}^{\infty}\log\frac{k_{j+1}}{k_j}(x)}$ defines a continuous, finite, positive definite metric $k_\infty$ on $\cO_{T_Y}(P_T)|_{g^{-1}(U\cap A^\circ)}$ which is the limit of $k_i$. The limiting weights are plurisubharmonic, so $k_\infty$ has semipositive curvature. By the triangle inequality and the inequality \eqref{eq:metric iteration contraction}, we deduce that for every $a\in U\cap A^\circ$,
\[
\begin{aligned}
 d_a\bigl(\Psi(k_\infty),k_\infty\bigr)&\le d_a\bigl(\Psi(k_\infty),\Psi(k_j)\bigr)+d_a(k_{j+1},k_\infty)\\
 &\le c\,d_a(k_\infty,k_j)+d_a(k_{j+1},k_\infty).
\end{aligned}
\]
Let $j\to\infty$. Then we deduce that $d_a\bigl(\Psi(k_\infty),k_\infty\bigr)=0$ and the two metrics agree on $T_{Y,a}$. Since it holds for every $a\in U\cap A^\circ$, we deduce that over $g^{-1}(U\cap A^\circ)$,
\[
 \Psi(k_\infty)=k_\infty.
\]

Over each fiber $T_{Y,a}$ where $a\in A^\circ$, the continuous, finite, positive definite metric on $\cO_{T_Y}(P_T)|_{T_{Y,a}}$ fixed by $\Psi$ is unique. Indeed, if $k$ and $l$ are two continuous, finite, positive definite metrics on $\cO_{T_Y}(P_T)|_{T_{Y,a}}$ fixed by $\Psi$, then
\[
 d_a(k,l)=d_a\bigl(\Psi(k),\Psi(l)\bigr)\le c\,d_a(k,l),
\]
and thus $d_a(k,l)=0$ and $k=l$. Thus, since the definition of $\Psi$ over each fiber is irrelevant to the choices of $U$ and $k_0$, the metrics $k_\infty$ which we just constructed locally glue to a continuous, finite, positive definite metric, still denoted by $k_\infty$, on $\cO_{T_Y}(P_T)|_{T_Y^\circ}$ where $T_Y^\circ:=g^{-1}(A^\circ)$.

\emph{Step 3.3.} We now extend the metric $k_\infty$ across $T_Y\setminus T_Y^\circ$. Keep working with the coordinate ball $U=B(0,5R)$ and the initial metric $k_0$ on $\cO_{T_Y}(P_T)|_{g^{-1}(U)}$. Let $U_1:=B(0,R)$ and $U_2:=B(0,2R)$. We next prove that there exists a constant $0<\varepsilon\le1$ such that over $g^{-1}(U_1\cap A^\circ)$,
\begin{equation}\label{eq:metric iteration bound}
 \Psi(k_0)\ge\varepsilon k_0.
\end{equation}
Write $dt:=dt_1\wedge\cdots\wedge dt_{\dim A}$ using the coordinates of the coordinate ball $U=B(0,5R)$. For $s\in H^0(U_2,\cE)$, let $\tilde s\in H^0(g^{-1}(U_2),\mathcal O_{T_Y}(D'_T))$ be the corresponding section under $\cE=g_*\cO_{T_Y}(D'_T)$. Let $\beta_s\in H^0(g^{-1}(U_2),\omega_{T_Y}\otimes L)$ be the image of $\tilde s\otimes (g^*dt)|_{g^{-1}(U_2)}$ under the natural isomorphism
\[
 \cO_{T_Y}(D'_T)\otimes g^*\omega_A\cong\omega_{T_Y}\otimes L.
\]
Recall that $L=\cO_{T_Y}(D'_T)\otimes\omega_{T_Y/A}^{-1}$. Here $g^*dt$ is a holomorphic frame of $g^*\omega_A$ over $g^{-1}(U)$.

We next prove by contradiction that there exists a constant $C>0$, independent of $s\in H^0(U_2,\cE)$, such that
\begin{equation}\label{eq:metric iteration section bound}
 \max_{x\in g^{-1}(\overline{U_1})}\left|\bigl(\rho(g^*s)\bigr)(x)\right|_{k_0,x}^2\le C\int_{g^{-1}(U_2)}|\beta_s|_{h_{B_T}k_0^c}^2.
\end{equation}
Here $|\beta_s|_{h_{B_T}k_0^c}^2$ is the nonnegative measurable $(\dim T_{Y},\dim T_{Y})$-form associated with $\beta_s$ and the metric $h_{B_T}k_0^c$, which is defined in \cite{HPS18}*{Section 14}. Suppose that no such constant $C>0$ exists. By rescaling sections, we can choose $s_j\in H^0(U_2,\cE)$ for every $j\ge1$ such that
\[
 \max_{x\in g^{-1}(\overline U_1)}\left|\bigl(\rho(g^*s_j)\bigr)(x)\right|_{k_0,x}^{2}=1,\qquad\int_{g^{-1}(U_2)}|\beta_{s_j}|_{h_{B_T}k_0^c}^{2}<\frac1j.
\]
Around any point in $g^{-1}(U_2)$, choose an open coordinate neighborhood contained in $g^{-1}(U_2)$ over which $L$ has a holomorphic frame $\ell$. Over the coordinate neighborhood, write
\[
 \beta_{s_j}=b_j\,dz_1\wedge\cdots\wedge dz_{\dim T_Y}\otimes\ell,\qquad |\ell|_{h_{B_T}k_0^c}^{2}=e^{-\psi}.
\]
Choose smaller open neighborhoods $W'$ and $W$ such that $\overline{W'}$ is a compact subset of $W$ and $\overline W$ is a compact subset of the chosen open coordinate neighborhood. Since $\psi$ is locally bounded from above, there exists a real number $M$ such that $\psi\le M$ on $W$. By the mean value inequality for $|b_j(x)|^2$ and the choice of $W'$ and $W$, we get a constant $C_{W',W}>0$, independent of $j$, such that
\[
 \begin{aligned}
 \max_{x\in\overline{W'}}|b_j(x)|^2&\le C_{W',W}\int_W|b_j|^2\,d\lambda\le C_{W',W}e^M\int_W|b_j|^2e^{-\psi}\,d\lambda\\
 &\le C_{W',W}e^M\int_{g^{-1}(U_2)}|\beta_{s_j}|_{h_{B_T}k_0^c}^{2}<\frac{C_{W',W}e^M}{j}.
 \end{aligned}
\]
Here $d\lambda$ denotes the Euclidean volume form in the chosen coordinates. We deduce that the local coefficients of $\beta_{s_j}$ converge uniformly to zero on compact subsets of $g^{-1}(U_2)$. Thus the same statement holds for $\tilde s_j$. Equivalently, $\tilde s_j\to 0$ in $H^0(g^{-1}(U_2),\cO_{T_Y}(D'_T))$ for the Fr\'echet topology, see \cite{GR09}*{Chapter VIII, Section A} and \cite{HPS18}*{Section 15} for related definitions and results regarding the Fr\'echet topology. By \cite{HPS18}*{Proposition 15.3}, $s_j\to 0$ in $H^0(U_2,\cE)$. Consider the map
\[
\begin{aligned}
    H^0(U_2,\cE)&\longrightarrow H^0(g^{-1}(U_2),\cO_{T_Y}(P_T)) \\
 s&\longmapsto\rho(g^*s).
\end{aligned}
\]
It is also the following map
\[
 H^0(U_2,\cE)\longrightarrow H^0(U_2,g_*\cO_{T_Y}(P_T))=H^0(g^{-1}(U_2),\cO_{T_Y}(P_T))
\]
induced by the morphism $\cE\to g_*\cO_{T_Y}(P_T)$ of coherent sheaves adjoint to $\rho$. The induced maps of morphisms of coherent sheaves on their section spaces are continuous with respect to the Fr\'echet topologies on the section spaces. Thus we deduce that the above map is continuous with respect to the Fr\'echet topologies by \cite{HPS18}*{Proposition 15.3}. Thus $\rho(g^*s_j)\to 0$ in $H^0(g^{-1}(U_2),\cO_{T_Y}(P_T))$. Equivalently, $\rho(g^*s_j)$ converges uniformly to zero on compact subsets of $g^{-1}(U_2)$. Since $g^{-1}(\overline U_1)$ is compact and $k_0$ is continuous, finite, and positive definite, we deduce that
\[
 1=\max_{x\in g^{-1}(\overline U_1)}\left|\bigl(\rho(g^*s_j)\bigr)(x)\right|_{k_0,x}^{2}\longrightarrow 0,
\]
which is a contradiction. Thus we have proved \eqref{eq:metric iteration section bound}. Note that \eqref{eq:metric iteration section bound} can also be proved by tedious estimates without invoking the Fr\'echet topologies. We use them here to make our arguments short and concise.

Fix $a\in U_1\cap A^\circ$ and $e\in\cE(a)$. Let $U':=B(a,3R)$. We have that $U_2\subseteq U'\subseteq U$. By applying the Ohsawa-Takegoshi extension theorem stated in \cite{HPS18}*{Theorem 14.4} to the morphism $g^{-1}(U')\to U'$ and the metric with semipositive curvature $h_{B_T}k_0^c$ on $L|_{g^{-1}(U')}$, we get a section $s\in H^0(U',\cE)$ with $s(a)=e$ and
\[
 \int_{g^{-1}(U')}|\beta_s|_{h_{B_T}k_0^c}^2\le\Vol(U')\|e\|_{H(k_0),a}^2,
\]
where $\beta_s$ is defined over $g^{-1}(U')$ by the same construction and $\Vol(U')$ denotes the volume of $U'$. For every $x\in T_{Y,a}$, we have $(\rho(g^*s))(x)=\rho_x(s(a))=\rho_x(e)$. By applying \eqref{eq:metric iteration section bound} to $s|_{U_2}$ and the above inequality, we deduce that
\[
 |\rho_x(e)|_{k_0,x}^2\le C\Vol (U')\|e\|_{H(k_0),a}^2.
\]
For a fixed $v\in\cO_{T_Y}(P_T)(x)$, take the minima of both sides of the above inequality over all $e$ such that $\rho_x(e)=v$. By the definition of the quotient metric $\Psi(k_0)$, we deduce that
\[
 |v|_{k_0,x}^2\le C\Vol (U')|v|_{\Psi(k_0),x}^2.
\]
Let $ \varepsilon:=\min\{1,\frac{1}{C\Vol (U')}\}$, and we have proved \eqref{eq:metric iteration bound}. By repeatedly applying $\Psi$ to $k_1\ge\varepsilon k_0$ over $g^{-1}(U_1\cap A^\circ)$, we deduce that $k_j\ge\varepsilon^{1+c+\cdots+c^{j-1}}k_0$ for every $j\ge1$. Let $j\to\infty$, and we deduce that $k_\infty\ge\varepsilon^{\frac{1}{1-c}}k_0$. Thus over $g^{-1}(U_1\cap A^\circ)$, the local weights of $k_\infty$ are bounded from above by the local weights of $k_0$ minus $\frac{\log\varepsilon}{1-c}$. Since $k_0$ is a continuous, finite, positive definite metric on $\cO_{T_Y}(P_T)|_{g^{-1}(U)}$, its weights are locally bounded from above over $g^{-1}(U_1)$ particularly. Thus the weights of $k_\infty$ are locally bounded from above near $g^{-1}(U_1)\setminus g^{-1}(U_1\cap A^\circ)$. The above arguments work for any point of $A$ and any sufficiently small coordinate ball centered around that point. Thus the weights of $k_\infty$ are locally bounded from above near $T_Y\setminus T_Y^\circ$. Since every irreducible component of $T_Y$ dominates $A$, $T_Y\setminus T_Y^\circ$ contains no irreducible component of $T_Y$. Thus by \cite{Dem12}*{Chapter I, Theorem 5.24}, the plurisubharmonic weights of $k_\infty$ extend uniquely across $T_Y\setminus T_Y^\circ$. The uniqueness guarantees that these extensions are compatible with each other. Thus we obtain a singular Hermitian metric $k_\infty$ with semipositive curvature on $\cO_{T_Y}(P_T)$ over the entire $T_Y$. It is continuous, finite, positive definite over $T_Y^\circ$.

\emph{Step 3.4.} We now prove that $\cQ$ admits a singular Hermitian metric with semipositive curvature and the minimal extension property and conclude our proof. The induced metric $h_{B_T}k_\infty^c$ on $L$ has semipositive curvature. Its multiplier ideal is trivial over $T_Y^\circ$ since $(T_Y,B_T)$ is klt and the weights of $k_\infty$ are locally bounded there. Thus the natural inclusion
\[
 g_*\bigl(\omega_{T_Y/A}\otimes L\otimes\mathcal I(h_{B_T}k_\infty^c)\bigr)\longrightarrow g_*(\omega_{T_Y/A}\otimes L)=\cE
\]
is an isomorphism over $A^\circ$. By \cite{HPS18}*{Theorem 21.1 and Corollary 21.2}, the sheaf $\cE$ admits a singular Hermitian metric with semipositive curvature and the minimal extension property. If $T_Y$ is disconnected, we apply our discussion to its connected components and take the orthogonal direct sum. Since $\pi:\cE\twoheadrightarrow\cQ$ is a surjection and $\cQ$ is torsion-free, the induced quotient metric on $\cQ$ has semipositive curvature and the minimal extension property by the same argument as at the beginning of our proof. Finally, every nonzero morphism $\cQ\to\cO_A$ splits by \cite{HPS18}*{Theorem 26.4}. The splitting is algebraic by GAGA.
\end{proof}

Next, we prove Theorem \ref{thm:gg for lc}. We will freely use the generic vanishing of $f_*\cO_X(D)$ and the description of $V_l^0(A,f_*\cO_X(D))$ where $f$ is a morphism from an lc pair $(X, \Delta)$ to an abelian variety $A$, $D\sim_{\Q}m(K_X+\Delta)$ is a Cartier divisor on $X$, and $m\ge1$ is a rational number. To be specific, under these assumptions, $f_*\cO_X(D)$ is a GV-sheaf and $V_l^0(A,f_*\cO_X(D))$ is a finite union of torsion subvarieties of $\Pic^0(A)$ for every $l\ge0$. When $m$ is an integer, they are well-known results by \cite{PS14b}*{Variant 5.5} and \cite{Shi16}*{Theorem 3.9} respectively. When $m$ is a rational number, they still hold and the proofs for them are almost verbatim. We omit the proofs for them here.

\begin{proof}[Proof of Theorem \ref{thm:gg for lc}]
By Lemma \ref{lem:modification}, we can assume that $(X, \Delta)$ is log smooth lc. We prove the statement by induction first on $\dim A$ and then on the number of lc centers of $(X, \Delta)$. The statement is immediate if $\dim A=0$. Now fix $\dim A$. Assume that the statement holds whenever the base abelian variety has dimension $<\dim A$. If the number of lc centers of $(X, \Delta)$ is $0$, then $(X, \Delta)$ is klt. Thus the statement follows from \cite{Men21}*{Theorem 1.1}. Next, we assume that the number of lc centers of $(X, \Delta)$ is $\ge1$.

Choose a minimal lc center $T$ of $(X, \Delta)$ and let $\cF:=f_*\cO_X(D)$, $\cK:=f_*(\cI_T\otimes\cO_X(D))$ and $\cQ:=\cF/\cK$. Let $\mu: Y:=\Bl_T X\to X$ be the blowup of $X$ along $T$ with exceptional divisor $E$. If $T$ is a divisor, let $\mu=\id_X$ and $E=T$. We have that $K_Y+\Delta_Y=\mu^*(K_X+\Delta)$. The coefficient of $E$ in $\Delta_Y$ is one, $\mu^*D-E\sim_{\Q}m(K_Y+\Delta_Y-E/m)$, and $(f\mu)_*\cO_Y(\mu^*D-E)=\cK$. Since $T$ is a minimal lc center, the lc centers of $(Y,\Delta_Y-E/m)$ are exactly the strict transforms of the lc centers of $(X, \Delta)$ other than $T$. By the induction hypothesis regarding the number of lc centers, $\cK$ is globally generated after pullback by an isogeny, and $\mathsf T(\cK)=\cK$. By Proposition \ref{prop:lifting from the lc locus}, the natural map
\[
 H^0(A,f_*\cO_X(D)\otimes\alpha)\longrightarrow H^0(A,\cQ\otimes\alpha)
\]
is surjective for every torsion line bundle $\alpha\in\Pic^0(A)$. Since $\cK\subseteq\cF$ and $\mathsf T(\cK)=\cK$, we deduce that $\cK=\mathsf T(\cK)\subseteq\mathsf T(\cF)$. By the above lifting property, we get the following exact sequence
\[
 0\longrightarrow\cK\longrightarrow\mathsf T(\cF)\longrightarrow\mathsf T(\cQ)\longrightarrow0.
\]
By Proposition \ref{prop:gg for lc restricted image m>1}, $\mathsf T(\cQ)$ a canonical direct summand. Thus $\mathsf T(\cQ)$ is a GV-sheaf. Since $\cK$ and $\cF$ are GV-sheaves, we deduce that $\mathsf T(\cF)$ and $\cU:=\cF/\mathsf T(\cF)$ are GV-sheaves.

Fix a fibration $p:A\to B$ to an abelian variety $B$ with $0\le\dim B<\dim A$, and a very ample line bundle $H$ on $B$. We will make choices of $p:A\to B$ to conclude our proof later. The sheaf $\cF\otimes p^*H$ is also a pushforward sheaf of a log smooth lc pair on $X$ under the morphism $f$. Thus $V_l^0(A,\cF\otimes p^*H)$ is a finite union of torsion subvarieties of $\Pic^0(A)$ for every $l\ge0$. By the induction hypothesis regarding the dimension of the base abelian variety, we deduce that for every torsion line bundle $\alpha\in\Pic^0(A)$, $\mathsf T(p_*(\cF\otimes\alpha))=p_*(\cF\otimes\alpha)$. Thus it is not hard to see that
\[
 p_*(\mathsf T(\cF)\otimes\alpha)=p_*(\cF\otimes\alpha),
\]
since $p:A\to B$ is a fibration between abelian varieties. Thus $H^0(A, \mathsf T(\cF)\otimes\alpha\otimes p^*H)=H^0(A, \cF\otimes\alpha\otimes p^*H)$ for every torsion line bundle $\alpha\in\Pic^0(A)$. We claim that the previous equality holds for every $\alpha\in\Pic^0(A)$. For every positive integer $l$, $V_l^0(A,\mathsf T(\cF)\otimes p^*H)\subseteq V_l^0(A,\cF\otimes p^*H)$. Since $V_l^0(A,\cF\otimes p^*H)$ is a finite union of torsion subvarieties of $\Pic^0(A)$, its torsion points are dense in itself. Since $V_l^0(A,\mathsf T(\cF)\otimes p^*H)$ is closed in $\Pic^0(A)$ and contains every torsion point of
$V_l^0(A,\cF\otimes p^*H)$, $V_l^0(A,\mathsf T(\cF)\otimes p^*H)=V_l^0(A,\cF\otimes p^*H)$ for every positive integer $l$. Thus for every $\alpha\in\Pic^0(A)$, we have that
\[
H^0(A, \mathsf T(\cF)\otimes\alpha\otimes p^*H)=H^0(A, \cF\otimes\alpha\otimes p^*H).
\]

Let $r$ be the minimum of $h^0(A, \mathsf T(\cF)\otimes\alpha\otimes p^*H)$ when varying $\alpha\in\Pic^0(A)$. We claim that
\[
 V^1(A,\mathsf T(\cF)\otimes p^*H)\subseteq V_{r+1}^0(A,\mathsf T(\cF)\otimes p^*H).
\]
Since $p^*H$ is nef, it is a GV-sheaf by \cite{Hac04}*{Theorem 1.2}. Thus $\mathsf T(\cF)\otimes p^*H$ is a GV-sheaf by \cite{PP11b}*{Theorem 3.2(1)}, since $\mathsf T(\cF)$ is also a GV-sheaf. Thus the Fourier--Mukai transform of $\mathsf T(\cF)\otimes p^*H$ is concentrated in degree $\dim A$ by \cite{PP11b}*{Theorem 2.3}. Denote it by $\cM[-\dim A]$, where $\cM$ is a coherent sheaf on $\Pic^0(A)$. Let $W:=\Pic^0(A)\setminus V_{r+1}^0(A,\mathsf T(\cF)\otimes p^*H)$ which is an open subset of $\Pic^0(A)$. For every $\alpha\in W$, $h^0(A, \mathsf T(\cF)\otimes\alpha\otimes p^*H)=r$. Let $W^{-1}:=\{\alpha^{-1}\mid\alpha\in W\}$. By the same argument as in \cite{LPS20}*{Proposition 2.3}, $\cM$ is locally free over $W^{-1}$, and for every $\alpha\in\Pic^0(A)$ and every $i\ge0$,
\[
 H^i(A,\mathsf T(\cF)\otimes\alpha\otimes p^*H)\cong\operatorname{Ext}^i_{\cO_{\Pic^0(A),\alpha^{-1}}}(\cM_{\alpha^{-1}},\C),
\]
where $\C$ is the residue field at $\alpha^{-1}\in\Pic^0(A)$. Since $\cM_{\alpha^{-1}}$ is free for every $\alpha\in W$, we deduce that $H^i(A,\mathsf T(\cF)\otimes\alpha\otimes p^*H)=0$ for every $\alpha\in W$ and every $i>0$. Thus
\[
V^1(A,\mathsf T(\cF)\otimes p^*H)\subseteq V_{r+1}^0(A,\mathsf T(\cF)\otimes p^*H)=V_{r+1}^0(A,\cF\otimes p^*H)\subsetneq\Pic^0(A).
\]
For every $\beta\in\Pic^0(B)$, $H\otimes\beta$ is globally generated. For every $\alpha\in\Pic^0(A)$ and every $\beta\in\Pic^0(B)$, consider the following exact sequence
\[
\begin{aligned}
0&\longrightarrow H^0(A,\mathsf T(\cF)\otimes\alpha\otimes p^*\beta\otimes p^*H)\longrightarrow H^0(A,\cF\otimes\alpha\otimes p^*\beta\otimes p^*H)\\
&\longrightarrow H^0(A,\cU\otimes\alpha\otimes p^*\beta\otimes p^*H)\longrightarrow H^1(A,\mathsf T(\cF)\otimes\alpha\otimes p^*\beta\otimes p^*H).
\end{aligned}
\]
From the above exact sequence, we deduce that
\begin{equation}\label{eq:inclusion}
\begin{aligned}
 V^0(A,\cU)+p^*\Pic^0(B)&\subseteq V^0(A,\cU\otimes p^*H)\subseteq V^1(A,\mathsf T(\cF)\otimes p^*H)\\
 &\subseteq V_{r+1}^0(A,\cF\otimes p^*H)\subsetneq\Pic^0(A).
\end{aligned}
\end{equation}

Recall that $V_l^0(A,\cF\otimes p^*H)$ is a finite union of torsion subvarieties of $\Pic^0(A)$ for every $l\ge0$. Let $Z$ be an irreducible component of $V^0(A,\cU)$. Consider the trivial fibration $p:A\to B$ where $B$ is a point. By \eqref{eq:inclusion}, $Z$ is contained in a torsion subvariety $C\subsetneq\Pic^0(A)$. We can assume that $C$ has the smallest dimension among all such torsion subvarieties containing $Z$. Write $C=\tau+C_0$ where $\tau\in\Pic^0(A)$ is a torsion point and $C_0\subsetneq\Pic^0(A)$ is an abelian subvariety. If $\dim C=\dim C_0>0$, by Poincar\'e's complete reducibility theorem, there exists an abelian subvariety $C_1\subsetneq\Pic^0(A)$ such that $C_0+C_1=\Pic^0(A)$, $C_0\cap C_1$ is finite, and the addition morphism $C_0\times C_1\to\Pic^0(A)$ is an isogeny. Let $B:=\Pic^0(C_1)$ and $p:A\to B$ the homomorphism dual to the inclusion $C_1\hookrightarrow\Pic^0(A)$. Then $p$ is a fibration such that $p^*\Pic^0(B)=C_1$. Since $\dim B=\dim A-\dim C_0$, we deduce that $0<\dim B<\dim A$. Moreover, we have that $C+p^*\Pic^0(B)=\tau+C_0+C_1=\Pic^0(A)$. By \eqref{eq:inclusion}, the irreducible closed set $Z+p^*\Pic^0(B)$ is contained in a proper torsion subvariety of $\Pic^0(A)$. This torsion subvariety cannot contain $C$, since otherwise it will contain $C$, $p^*\Pic^0(B)$ and thus $\Pic^0(A)$ which is a contradiction. The intersection of this torsion subvariety and $C$ is a finite union of torsion subvarieties with dimension smaller than $\dim C$. One of these smaller torsion subvarieties contains $Z$, which contradicts the choice of $C$. Thus $V^0(A,\cU)$ is a finite set consisting of torsion points.

If $\cU\ne0$, we deduce that $V^0(A,\cU)$ is not empty by \cite{HPS18}*{Lemma 7.4} since $\cU$ is a GV-sheaf. Choose $\alpha^{-1}\in V^0(A,\cU)$. Since $V^0(A,\cU)$ is a finite set consisting of torsion points, $\alpha^{-1}$ is an isolated torsion point. By \cite{PP11a}*{Proposition 3.15}, $\alpha^{-1}\in V^{\dim A}(A,\cU)$. By Serre duality, $\Hom_A(\cU,\alpha)\cong H^{\dim A}(A,\cU\otimes\alpha^{-1})^\vee\ne0$. Thus we obtain a nonzero morphism $\cU\to\alpha$ to the torsion line bundle $\alpha$. By the exact sequences $0\to\cK\to\cF\to\cQ\to0$ and $0\to\cK\to\mathsf T(\cF)\to\mathsf T(\cQ)\to0$, we deduce that $\cU=\cF/\mathsf T(\cF)\cong\cQ/\mathsf T(\cQ)$. Thus we obtain a nonzero morphism $\cQ\to\cQ/\mathsf T(\cQ)\to\alpha$ to the torsion line bundle $\alpha$. Then we deduce that $\Supp\cQ=A$. By the definition of $\cQ$, $\Supp\cQ\subseteq f(T)$. Thus $T$ dominates $A$. By applying Theorem \ref{thm:image splitting} to $D-f^*\alpha\sim_{\Q}m(K_X+\Delta)$, we deduce that the nonzero morphism $\cQ\to\cQ/\mathsf T(\cQ)\to\alpha$ splits. However, the image of the morphism $\alpha\to\cQ$ which splits $\cQ\to\cQ/\mathsf T(\cQ)\to\alpha$ is contained in $\mathsf T(\cQ)$. Thus the image of $\alpha\to\cQ\to\cQ/\mathsf T(\cQ)\to\alpha$ is $0$, which is a contradiction. Thus $\cU=0$ and $\cF=\mathsf T(\cF)$. Thus $f_*\cO_X(D)$ is globally generated after pullback by an isogeny.
\end{proof}

Finally, let us prove the following direct corollaries of Theorem \ref{thm:gg for lc}.

\begin{coro}\label{coro:existence of singular Hermitian metrics with semipositive curvature and the minimal extension property}
Let $f$ be a surjective morphism from an lc pair $(X, \Delta)$ to an abelian variety $A$, $m>1$ a rational number and $D$ a Cartier divisor on $X$ such that $D\sim_{\Q}m(K_X+\Delta)$. Then $f_*\mathcal{O}_X(D)$ admits a singular Hermitian metric with semipositive curvature and the minimal extension property.
\end{coro}

\begin{proof}
Since $f$ is surjective, $f_*\mathcal{O}_X(D)$ is torsion-free. By Theorem \ref{thm:gg for lc}, $\mathsf T(f_*\mathcal{O}_X(D))=f_*\mathcal{O}_X(D)$, and $f_*\mathcal{O}_X(D)$ is a torsion-free quotient sheaf of a finite direct sum of torsion line bundles on $A$ which admits a singular Hermitian metric with semipositive curvature and the minimal extension property. Thus the induced quotient metric on $f_*\mathcal{O}_X(D)$ satisfies the same properties.
\end{proof}

\begin{coro}\label{coro:real coefficient gg for lc}
Conjecture \ref{conj:global generation} holds when $m>1$ is a real number, $(X, \Delta)$ is an lc pair with real coefficients, and $D$ is a Cartier divisor on $X$ such that $D\sim_{\R}m(K_X+\Delta)$.
\end{coro}

\begin{proof}
By Lemma \ref{lem:modification}, we can assume that $(X, \Delta)$ is log smooth lc. Since $D$ is Cartier, there exist a rational number $m'>1$ and an lc pair $(X,\Delta')$ with rational coefficients such that $D\sim_{\Q}m'(K_X+\Delta')$ by the standard rational approximation techniques in the minimal model program, see \cite{Kol23}*{Section 11.4} for details. Thus $f_*\cO_X(D)$ is globally generated after pullback by an isogeny by Theorem \ref{thm:gg for lc}.
\end{proof}

\begin{rem}
Many results and arguments in this section still work when the base variety is an arbitrary smooth (not necessarily quasi-projective) variety and not necessarily an abelian variety by the same methods as in this section. We will discuss them in our future work.
\end{rem}

\section{Counterexamples}\label{sec:Counterexamples}

In this short section, we provide the promised counterexamples. This section can be read independently of Section \ref{sec:Positivity of pushforwards of lc pairs}.

\begin{proof}[Proof of Theorem \ref{thm:counterexample m>=1}]
Fix a rational number $m\ge1$, and write $m=a/b$, where $a\ge b\ge1$ are coprime integers. Let $E$ be an elliptic curve with origin $o$, and choose different torsion points $t_1,\ldots,t_b\in E$, with $t_1=o$. Let
\[
 A:=E\times E,\qquad T:=E\times\{o\},\qquad C:=\{o\}\times E,\qquad T_j:=E\times\{t_j\},\qquad S:=\sum_{j=1}^bT_j.
\]
Since the points $t_j$ are torsion, $S\sim_{\Q}bT$. The scheme theoretic intersection $Z:=S\cap C$ consists of the $b$ different reduced points $(o,t_1),\ldots,(o,t_b)$.

Let $\pi:X\to A$ be the blowup of $A$ along $Z$, and let $R_j$ be the exceptional curve over $(o,t_j)$. We use primes to denote strict transforms and let
\[
 \Delta:=S'+\frac ba C',\qquad D:=\pi^*(aT+C)-\sum_{j=1}^bR_j.
\]
The irreducible components of $\Delta$ are smooth and pairwise disjoint, so $(X,\Delta)$ is log smooth plt. Moreover,
\[
 K_X+\Delta\sim_{\Q}\pi^*\left(S+\frac ba C\right)-\frac ba\sum_{j=1}^bR_j\sim_{\Q}\frac ba\left(\pi^*(aT+C)-\sum_{j=1}^bR_j\right)=\frac ba D.
\]
Thus $D$ is a Cartier divisor on $X$ and $D\sim_{\Q}m(K_X+\Delta)$. Consider the external tensor product $L:=\cO_A(aT+C)\cong\cO_E(o)\boxtimes\cO_E(ao)$. We have that $\pi_*\cO_X(\ell D)\cong\cI_Z^\ell\otimes L^{\otimes\ell}$ for every integer $\ell\ge1$.

By direct calculations, $D$ is nef and big. Thus $K_X+\Delta$ is nef, big, and thus semiample by the abundance theorem in dimension $2$. Moreover, $2_A^*\pi_*\cO_X(\ell D)$ is globally generated for every $\ell\ge1$, where $2_A:A\to A$ is multiplication by $2$ on $A$. By the definition of $Z$, we have the following surjection
\[
 \bigl(L\otimes\cO_A(-S)\bigr)\oplus\bigl(L\otimes\cO_A(-C)\bigr)\twoheadrightarrow\pi_*\cO_X(D).
\]
After pullback by $2_A$, the two line bundles on the left hand side are external tensor products of line bundles of degrees $(4,4(a-b))$ and $(0,4a)$. Moreover, their degree zero factors are trivial by our choices of these divisors (note that it can happen that $a-b=0$). Both line bundles are globally generated. Thus $2_A^*\pi_*\cO_X(D)$ is globally generated. Since the multiplication $\Sym^{\ell}(\pi_*\cO_X(D))\twoheadrightarrow\cI_Z^\ell\otimes L^{\otimes\ell}$ is surjective for every $\ell\ge1$, $2_A^*\pi_*\cO_X(\ell D)$ is globally generated for every $\ell\ge1$. In particular, Conjecture \ref{conj:global generation} holds for $\pi_*\cO_X(D)$. Thus for every $\ell\ge1$, $\mathsf T(\pi_*\cO_X(\ell D))=\pi_*\cO_X(\ell D)$, and $\pi_*\cO_X(\ell D)$ is a torsion-free quotient sheaf of a finite direct sum of torsion line bundles on $A$ which admits a singular Hermitian metric with semipositive curvature and the minimal extension property. Thus the induced quotient metric on $\pi_*\cO_X(\ell D)$ satisfies the same properties for every $\ell\ge1$.

Next, we show that $\pi_*\cO_X(D)$ does not admit any Chen--Jiang decomposition, even after pullback by an isogeny. For every $\alpha\in\Pic^0(A)$, the line bundle $L\otimes\alpha$ is ample. By Kodaira vanishing, we deduce that
\[
H^1\bigl(A,\pi_*\cO_X(D)\otimes\alpha\bigr)\cong\Coker\left(H^0(A,L\otimes\alpha)\longrightarrow H^0\bigl(Z,(L\otimes\alpha)|_Z\bigr)\right).
\]
If $\alpha=\cO_E\boxtimes\beta$ with $\beta\in\Pic^0(E)$, every section of $L\otimes\alpha$ vanishes on $Z$, since every section of $\cO_E(o)$ vanishes at $o$. Thus
\[
\{\cO_E\}\times\Pic^0(E)\subseteq V^1\bigl(A,\pi_*\cO_X(D)\bigr).
\]
For any isogeny $\varphi:A'\to A$, let $\hat\varphi$ be the dual isogeny of $\varphi$. By \eqref{eq:isogeny}, we deduce that $V^i(A',\varphi^*\pi_*\cO_X(D))=\hat\varphi\bigl(V^i(A,\pi_*\cO_X(D))\bigr)$. Thus $\codim_{\Pic^0 (A')} V^1(A',\varphi^*\pi_*\cO_X(D))\le1$, so $\varphi^*\pi_*\cO_X(D)$ is not M-regular. Since $\varphi$ is \'etale,
\[
 \varphi^*\pi_*\cO_X(D)\cong\cI_{\varphi^{-1}(Z)}\otimes\varphi^*L
\]
is torsion-free of rank one and is not locally free. Suppose that it admits a Chen--Jiang decomposition. Every nonzero direct summand is torsion-free of positive rank, so the decomposition can have only one direct summand. Thus
\[
 \varphi^*\pi_*\cO_X(D)\cong\alpha_1\otimes p_1^*\cF_1,
\]
where $A_1$ is an abelian variety, $p_1:A'\to A_1$ is a fibration, $\cF_1$ is a nonzero M-regular coherent sheaf on $A_1$, and $\alpha_1\in\Pic^0(A')$ is a torsion line bundle. Since $p_1$ is flat and surjective, $\cF_1$ is torsion-free of rank one. If $\dim A_1\le1$, then $\cF_1$ is locally free, contradicting the fact that $\varphi^*\pi_*\cO_X(D)$ is not locally free. If $\dim A_1=2$, then $p_1$ is an isomorphism, and $\varphi^*\pi_*\cO_X(D)$ would be M-regular, again a contradiction. Thus $\varphi^*\pi_*\cO_X(D)$ does not admit any Chen--Jiang decomposition.

In particular, $\pi_*\cO_X(D)$ does not admit any Chen--Jiang decomposition. If $\pi_*\cO_X(D)$ were isomorphic to a direct summand of $h_*\omega_V$ for a projective morphism $h:V\to A$ from a smooth variety $V$, then it would admit such a decomposition by \cite{PPS17}*{Theorem A} and \cite{LPS20}*{Proposition 3.6}. Thus $\pi_*\cO_X(D)$ is not a canonical direct summand either. Thus Conjectures \ref{conj:CJ decomposition} and \ref{conj:canonical direct summand} both fail for any rational number $m\ge1$.
\end{proof}

\begin{proof}[Proof of Theorem \ref{thm:counterexample m=1}]
Take $a=b=1$ in the construction in the proof of Theorem \ref{thm:counterexample m>=1}. Thus $\pi:X\to A=E\times E$ is the blowup of $p=(o,o)$, with exceptional curve $R$, and
\[
 \Delta=T'+C',\qquad D=\pi^*(T+C)-R\sim K_X+\Delta.
\]
The pair $(X,\Delta)$ is log smooth plt, and $K_X+\Delta$ is nef, big, and semiample.

Let $s:A\to E$ be the addition morphism and $f:=s\pi$. Then $f_*\cO_X\cong s_*\cO_A\cong\cO_E$, and $f$ is a fibration. We have that
\[
 f_*\cO_X(D)\cong s_*\bigl(\cI_p\otimes\cO_A(T+C)\bigr).
\]
By twisting the Koszul complex of $p=T\cap C$ by $\cO_A(T+C)$, we have that
\[ 0\longrightarrow\cO_A\longrightarrow\cO_A(T)\oplus\cO_A(C)\longrightarrow\cI_p\otimes\cO_A(T+C)\longrightarrow0.
\]
The addition morphism becomes a product projection after being composed with an automorphism of $A$, so $s_*\cO_A\cong R^1s_*\cO_A\cong\cO_E$. Moreover, each of $\cO_A(T)$ and $\cO_A(C)$ restricts to a degree one line bundle on every fiber of $s$. By cohomology and base change, we have that
\[
 s_*\cO_A(T)\cong s_*\cO_A(C)\cong\cO_E,\qquad R^1s_*\cO_A(T)=R^1s_*\cO_A(C)=0.
\]
Consider the pushforward of the above Koszul complex under the morphism $s$. By the above isomorphisms derived from cohomology and base change, $s_*\cO_A\to s_*\cO_A(T)\oplus s_*\cO_A(C)$ is identified with the map $\cO_E\to\cO_E\oplus\cO_E$ given, up to signs, by the constant vector $(1,-1)$ whose cokernel is $\cO_E$. Thus the pushforward sequence gives
\[
0\longrightarrow\cO_E\longrightarrow f_*\cO_X(D)\longrightarrow\cO_E\longrightarrow0.
\]
Every section of $\cO_A(T+C)\cong\cO_E(o)\boxtimes\cO_E(o)$ vanishes at $p$, thus
\[
h^0\bigl(E,f_*\cO_X(D)\bigr)=h^0\bigl(A,\cI_p\otimes\cO_A(T+C)\bigr)=1.
\]
Thus the above short exact sequence is non-split. Thus $f_*\cO_X(D)$ is the unique non-split extension of $\cO_E$ by $\cO_E$. By \cite{Men23}*{Example 4.3}, $f_*\cO_X(D)$ does not admit any Chen--Jiang decomposition and is not globally generated after pullback by any isogeny. Thus Conjecture \ref{conj:global generation} fails for $m=1$.

Finally, suppose that $f_*\cO_X(D)$ admits a singular Hermitian metric with semipositive curvature and the minimal extension property. By \cite{HPS18}*{Theorem 26.4} and GAGA, the nonzero quotient morphism $f_*\cO_X(D)\to\cO_E$ in the above short exact sequence splits, a contradiction. Thus $f_*\cO_X(D)$ does not admit any such metric.
\end{proof}

\bibliographystyle{amsalpha}
\bibliography{biblio}
\end{document}